\documentclass[12pt]{article}
\usepackage{epsfig}
\usepackage{latexsym}
\usepackage{amsfonts}
\usepackage{amsmath}
\usepackage[psamsfonts]{amssymb}

\def\R{\mathbf{R}}
\def\bC{\mathbf{\overline{C}}}

\def\Z{\mathbf{Z}}
\def\N{\mathbf{N}}

\newtheorem{thm}{Theorem}[section]
\newtheorem{example}[thm]{Example}
\newtheorem{prop}[thm]{Proposition}
\newtheorem{lemma}[thm]{Lemma}

\newtheorem{cor}[thm]{Corollary}
\newtheorem{rmk}[thm]{Remark}
\newtheorem{df}[thm]{Definition} 

\numberwithin{equation}{section}

\begin{document}
\title{Spherical rectangles}
\author{Alexandre Eremenko\thanks{Supported by NSF grant DMS-1361836.}
\
and Andrei Gabrielov\thanks{Supported by NSF grant DMS-1161629.}}
\maketitle
\begin{abstract}
We study spherical quadrilaterals whose angles are
odd multiples of $\pi/2$, and the equivalent accessory parameter problem
for the Heun equation.
We obtain a classification of these quadrilaterals up to isometry. For given
angles, there are finitely many one-dimensional
continuous families which we enumerate.
In each family the conformal modulus is either bounded from above
or bounded from below, but not both, and the numbers of families
of these two types are equal.
The results can be translated
to classification of Heun's equations with real parameters, whose
exponent
differences are odd multiples of $1/2$, with unitary monodromy.


\end{abstract}

\section{Introduction}
A (marked) {\em circular polygon} is a closed disk $Q$ with marked boundary points
$a_0,\ldots,a_{n-1}$, which are called {\em corners}, enumerated cyclically according to the positive
orientation of $\partial Q$,
equipped with a Riemannian metric of constant curvature
$1$ with conic singularities at the corners, and such that each {\em side}
$(a_j,a_{j+1})$ has constant geodesic curvature. Two such polygons $Q$
and $Q'$ are {\em congruent}
if there is an orientation-preserving isometry
between them which sends each corner $a_j$ of $Q$ to
the corner $a_j^\prime$ of $Q'$.

Polygons with $n=2,\ 3$ and $4$ are called {\em digons, triangles} and {\em quadrilaterals}, respectively.

If each side has zero geodesic curvature then $Q$ is
called a {\em spherical polygon}.
At every corner $a_j$, an interior angle $\alpha_j\geq 0$
is defined, and in what follows
we measure all angles in half-turns.
So angle $\alpha$ means an angle of $\pi\alpha$ radians,
in particular, ``integer angle'' is an integer multiple of $\pi$ radians.
A circular quadrilateral $Q$ whose angles are odd multiples of $1/2$
is called a {\em circular rectangle}. If such a quadrilateral is spherical,
it is called a {\em spherical rectangle}. In this paper we describe
the set of spherical rectangles with prescribed angles.

As every surface of positive curvature $1$ is locally isometric to a piece
of the unit sphere, every circular polygon can be described in terms
of the {\em developing map} $f:Q\to\bC$ which is an analytic function
on $Q\setminus\{ a_0,\ldots,a_{n-1}\}$ mapping every side into
a circle on the Riemann sphere. For spherical polygons these circles
are geodesic (great circles). This function $f$ is a local homeomorphism
at each point except the corners, and at a corner $a$ satisfies
$$f(z)-f(a)\sim c(z-a)^\alpha,$$
where $\alpha>0$ is the angle at this corner. (If $\alpha=0$,
the right-hand side has to be replaced by $c/\log(z-a)$.)

Each such function defines a circular polygon by the pull-back of the spherical
metric from $\bC$ to $Q$. If none of the $\alpha_j$ equals $1$,
then the pair $(f,a_0)$ defines the polygon uniquely.
Two pairs $(f_1,a_0)$ and $(f_2,a_0^\prime)$ define
congruent polygons if $f_2=\psi\circ f_1\circ\phi$,
where $\psi$ is a rotation of the Riemann sphere, and $\phi$
a conformal automorphism of the disk with the property $\phi(a_j)=a_j^\prime$ for all $j$.

This paper is a part of the project whose goal is to understand metrics
of constant positive curvature with conic singularities on compact surfaces.
\cite{Troy2,LT,E,EGT1,EGT2,EGT3,MP,Lin1,Lin2}. An important class of
such metrics can be obtained by gluing a spherical polygon to its
mirror image isometrically along the boundary. Metrics of positive curvature
on the sphere obtained in this way are characterized by the symmetry
property: all conic singularities belong to a circle on the sphere,
and the metric is symmetric with respect to this circle.

In this paper, we classify spherical rectangles.
The cases when at least one of the angles of a spherical quadrilateral
is integer were considered in \cite{EGT3,EGT2,EGT1} and \cite{EGSV}.

If we use the upper half-plane as $Q$, then
the developing map of every circular quadrilateral is a ratio of two linearly
independent solutions of the Heun equation
\begin{equation}\label{heun}
y''+\left(\sum_{j=0}^2\frac{1-\alpha_j}{z-a_j}\right)y'+
\frac{\alpha'\alpha'' z-\lambda}{\prod_{j=0}^2(z-a_j)}y=0,
\end{equation}
where
\begin{eqnarray}\label{heun2}
\alpha'&=&(2+\alpha_3-\alpha_2-\alpha_1-\alpha_0)/2,\\
\alpha''&=&(2-\alpha_3-\alpha_2-\alpha_1-\alpha_0)/2,\label{heun3}
\end{eqnarray}
with the standard normalization $(a_0,a_1,a_2,a_3)=(0,1,a,\infty)$, $a\in(1,+\infty)$.
Here $\alpha_j$ are the angles at the corners, and $\lambda$ is a real
{\em accessory parameter}. Each pair of linearly independent solutions
of such an equation defines a circular quadrilateral. Different pairs of
linearly independent solutions of the same equation define {\em equivalent}
quadrilaterals: their developing maps are related by post-composition with
linear-fractional transformations. We will
later see that an equivalence class
may contain at most one spherical quadrilateral, up to congruence.

The condition that an equivalence class contains a spherical quadrilateral
translates to the following condition on the Heun equation: the projective
monodromy group must be conjugate to a subgroup of $SU(2)$. So the problem
of classification of spherical quadrilaterals with prescribed corners and
angles is equivalent to the problem of classification of Heun's equation
with prescribed $a_j$ and $\alpha_j$ whose monodromy is {\em unitarizable},
that is conjugate to a subgroup of $SU(2)$. The correspondence
between the metrics on the sphere and Heun's equations with unitarizable
monodromy is bijective. Each symmetric metric on the sphere corresponds
to two spherical quadrilaterals which are related by an anti-conformal involution,
and to a unique normalized Heun equation with real $a,\lambda$ and unitarizable
monodromy.

In the case that all angles are odd multiples of $1/2$, Heun's equation
can be explicitly solved in terms of elliptic integrals.
This fact was discovered by Darboux \cite{Darboux} who generalized
Hermite's work \cite{Hermite} for the Lam\'e equation. For the study of
general circular rectangles in connection with Heun's equation with
$\alpha_j$ multiples of $1/2$ we refer to the paper of Van Vleck \cite{VV}.

We recall how this explicit solution is obtained.
\vspace{.1in}

\noindent
{\bf Theorem A.} {\em Suppose that all $\alpha_j$ in
(\ref{heun})--(\ref{heun3})
are odd multiples of $1/2$. Then there are two linearly independent
solutions of (\ref{heun}) whose ratio is of the form
\begin{equation}\label{solution}
f(z)=\exp(I(z)),
\end{equation}
where 
\begin{equation}\label{integral}
I(z)=\int_{z_0}^z\prod_{j=0}^3(\zeta-a_j)^{\alpha_j-1}
\frac{d\zeta}{P(\zeta)}
\end{equation}
and $P(z)=P(z;a,\lambda)$ is a real polynomial in all three variables.
This polynomial satisfies
the third order linear differential equation
\begin{equation}\label{hermite}
w'''+3pw''+(p'+2p^2+4q)w'+(4pq+2q')w=0,
\end{equation}
where $p$ and $q$ are the coefficients in front of $y'$ and $y$ in
(\ref{heun}).}
\vspace{.1in}

Equation (\ref{hermite}) has one-dimensional space of polynomial solutions
of degree
\begin{equation}\label{degree}
\deg P=\sum_{j=0}^3\alpha_j-2.
\end{equation}
This permits to find $P$ by rational operations.
That $P$ satisfies (\ref{hermite}) guarantees that all residues of 
the integrand in (\ref{integral}) are of the form $\pm c$ with some real $c$.
The condition $c=1$ defines $P$ up to a sign.

Periods of the integral (\ref{integral}), other than those coming from the
residues, form a lattice generated by two canonical periods (integrals over
adjacent real segments).
Of these two canonical periods one is real and another is pure imaginary.
The condition that the monodromy of (\ref{heun}) is unitarizable
means that both periods must be imaginary, therefore the
real period must vanish.
For a fixed real $a$, and given angles, this imposes a transcendental
equation on $\lambda$. It is not clear how to determine or estimate
the number of real solutions of this equation, but for small
angles $\alpha_k$ it can be solved numerically. The results of
computation are described in Section \ref{limits}.

Instead we use a geometric method which allows us to classify
spherical rectangles, and describe their geometry.
The following elementary statement was proved in \cite{EG2}.

\begin{prop}\label{onecircle} Let $f$ be the developing map of a
spherical rectangle $Q$.
Then there are two opposite sides of $Q$ whose $f$-images are
contained in the same circle, and the other pair of opposite sides
is mapped to distinct circles.
\end{prop}

Thus the boundary of a spherical rectangle $Q$ is mapped by $f$ to the union
of three great circles, one of them, say $C$, being orthogonal to the
other two, $C'$ and $C''$.
Let $\theta\in(0,1)$ be the angle between the circles $C'$ and $C''$.
There are two choices for this angle (the other one being $1-\theta$).
See Definition \ref{df:theta} below for the unique choice of the angle $\theta$
associated with a spherical rectangle $Q$.

The $f$-preimage in $Q$ of the three circles is called the {\em net} of $Q$.
The net is a combinatorial invariant of $Q$,
defined up to an orientation-preserving
homeomorphism of $Q$ respecting its initial corner $a_0$.

\begin{prop}\label{isometry}
A marked spherical rectangle $Q$ is defined uniquely up to
isometry by its net and the angle $\theta\in(0,1)$.
\end{prop}
\vspace{.1in}

This will be proved in Section \ref{nets}. In section \ref{nets}, we will explicitly
describe all possible nets of spherical rectangles (Theorems \ref{classification} and \ref{existence}).
As a consequence we will obtain the following necessary and sufficient conditions on the angles of
a spherical rectangle:

\noindent
\begin{thm}\label{theorem1}
Let $A_0,\ldots,A_3$ be non-negative integers, and
\begin{equation}\label{delta}
\delta=(A_1+A_3-A_0-A_2)/2.
\end{equation}
Then for the existence of a spherical rectangle with angles $\alpha_j=A_j+1/2$
it is necessary and sufficient that one of the following conditions
be satisfied: either
\begin{equation}\label{1}
\delta\geq 1,\quad A_1\geq 1,\quad A_3\geq 1,
\end{equation}
or
\begin{equation}\label{2}
\delta\leq -1,\quad A_0\geq 1,\quad A_2\geq 1.
\end{equation}
\end{thm}
\vspace{.1in}

To state our next result we need some definitions.
We recall that we consider {\em marked} spherical rectangles.
Two of them are congruent if there is an orientation-preserving
isometry of $Q$ respecting the initial corner $a_0$.

Every quadrilateral can be mapped conformally onto a flat rectangle
with vertices $(0,1,1+iK,iK)$ and all angles $1/2$, such that  $a_0$ maps
to $0$. The number $K$ is called the {\em modulus} of the quadrilateral.

Each pair $(\Gamma,\theta)$, where $\Gamma$ is a net and $\theta\in(0,1)$, defines a marked spherical rectangle
$Q(\Gamma,\theta)$ (see Theorem \ref{existence}). Thus the set of all spherical rectangles
with given angles consists of curves $\theta\mapsto Q(\Gamma,\theta)$
parameterized by $\theta$ and labeled by the nets. The modulus $K$
of $Q(\Gamma, \theta)$ is a continuous function of $\theta$.
There are two kinds of these curves:

On the curves of the first kind, $K\to 0$ as $\theta\to 0$,
while $K$ tends to a non-zero limit $K_{\mathrm{crit}}(\Gamma)$
as $\theta\to 1$.

On the curves of the second kind, $K\to+\infty$ as $\theta\to 0$,
while $K$ tends to a non-zero limit $K_{\mathrm{crit}}(\Gamma)$ as $\theta\to 1$.

This is proved in Theorem \ref{limit}, and we give few examples of computation
of the limits $K_{\mathrm{crit}}$ in Section \ref{limits}.
In all our examples $K$ is a monotone
function of $\theta$. This is proved in \cite{EG2} for the simplest family
of spherical quadrilaterals with angles $(3/2,1/2,3/2,1/2)$ but
it is unlikely that this property holds in general.
However, it is true for sufficiently small and large values of $K$.

\begin{prop}\label{small-large}
Each curve $Q(\Gamma, \theta)$ has finitely many intervals on which $K$ is monotone.
In particular, for sufficiently small (resp., large) $K>0$, there is a unique spherical
rectangle  with the modulus $K$ in a curve $Q(\Gamma, \theta)$ of the first (resp., second) kind.
\end{prop}

This follows from the general theory of o-minimal structures (see, e.g., \cite{vdd}), since
the integral in (\ref{integral}) is a Pfaffian function in the sense of Khovanskii \cite{khov}.
It was shown in \cite{sp} that the structure generated by Pfaffian functions is o-minimal.

Our final results count the nets for spherical rectangles with given angles.

The quadruple $(A_0,\dots,A_3)$ is {\em special} if $\delta$ in (\ref{delta})
is an odd integer and one of the following holds:
either $A_1\ge\delta>0$ and $A_3\ge\delta>0$ or $A_0\ge-\delta>0$ and $A_2\ge-\delta>0$.
\vspace{.1in}

Define
\begin{equation}\label{K1}
M_1=\left[\min\left\{\frac{A_1+1}2,\frac{A_3+1}2,\frac{1+\delta}2\right\}\right],
\end{equation}
\begin{equation}\label{K2}
M_2=\left[\min\left\{\frac{A_0+1}2,\frac{A_2+1}2,\frac{1-\delta}2\right\}\right],
\end{equation}
\begin{equation}\label{N}
N=\min\left\{A_0+\frac{1+\delta}{2},A_1+\frac{1-\delta}{2},
A_2+\frac{1+\delta}{2},A_3+\frac{1-\delta}{2}\right\}.
\end{equation}
Note that conditions (\ref{1}) and (\ref{2}) are satisfied when $M_1>0$ and $M_2>0$, respectively.

\noindent
\begin{thm}\label{theorem2}
For a special quadruple $(A_0,\ldots,A_3)$ satisfying either (\ref{1}) or (\ref{2})
there exist $2N$ one-parametric families of congruence classes of
marked spherical rectangles with angles $A_j+1/2$.
If $(A_0,\ldots,A_3)$ is not special but satisfies (\ref{1}) (resp., (\ref{2})) then
there exist $2M_1$ (resp., $2M_2$) one-parametric families of congruence classes of
marked spherical rectangles with angles $A_j+1/2$.

Each family is parameterized by $\theta\in(0,1)$
(see Definition \ref{df:theta}).
Each family contains either rectangles of arbitrarily small moduli or
arbitrarily large moduli but not both. The numbers of families of
both types are equal, so for each type this number is either $N$ or
$M_1$ or $M_2$, depending on $(A_0,\dots,A_3)$.
\end{thm}
\vspace{.1in}

\begin{rmk}\label{exact}
Theorem \ref{theorem2} and Proposition \ref{small-large} imply that
the number of spherical rectangles with given angles $A_j+1/2$
is exactly $N$ or $M_1$ or $M_2$, depending on $(A_0,\dots,A_3)$,
for sufficiently small and large values of $K$.
\end{rmk}

For fixed angles $\alpha_j\in \N^++1/2$,
consider the two-parametric family of Heun's
equations (\ref{heun}) with parameters
$(a,\lambda)\in(1,+\infty)\times\R.$
Equivalence classes of circular rectangles are in
correspondence to such Heun's equations.
One-parametric families of spherical rectangles of
Theorem \ref{theorem2} correspond to smooth disjoint curves
in the half-plane $(a,\lambda)$, each having one end in this half-plane.
On the other end, either $a\to 1$ or $a\to\infty$.
\vspace{.1in}

When $\theta=p/q$ is rational, the monodromy group
of the developing map is finite,
so $f$ is an algebraic function. In this case, $f=g^{-1}\circ h$, where
$$g=-\frac14\left(z^q+\frac1{z^q}-2\right),$$
and $h$ is a rational Belyi function,
which means that the only critical values of $h$ are
$0,1,\infty$, \cite{Schneps}. Function $g$ is also a Belyi function,
it is called the fundamental rational function of the dihedral group
\cite{Klein2}. The set $g^{-1}(\R)$ consists of the unit circle and $q$
lines $\{ z=t\exp(\pi ik/q): t\in\R,\;k=0,\ldots,q-1\}.$
In the simplest case $\theta=1/2$, the image $f(\partial Q)$ is
contained in the union of the unit circle $C$, real line $C'$, and imaginary line $C''$.
The monodromy is the Klein Viergroup
$\Z^2\times \Z^2$, represented as $\{ z,-z,1/z,-1/z\}$ and $g$
 has the property that $g^{-1}(\R)=C\cup C'\cup C''$.
Then it is easy to see that our net, together with its
reflection in the real line, coincides with
$h^{-1}(\R)$. The set $h^{-1}(\R)$ for a Belyi function $h$
is a triangulation of the sphere with all vertices of even degree.
Therefore our classification of the nets can be restated as classification
of triangulations $\mathbf T$ of the sphere with the following properties:

a) $\mathbf T$ is symmetric with respect to $\R$, and $\R$ is contained in the $1$-skeleton of $\mathbf T$,

b) There are four vertices $a_j$ of $\mathbf T$ on the real line of prescribed orders $4A_j+2$.

c) All other vertices of $\mathbf T$ have order $4$.

For $q\in\{2,3\}$ and $(A_0,\ldots,A_3)=(1,0,1,0)$, algebraic developing maps are explicitly written in \cite{EG2}.

\noindent
\section{Spherical polygons with the sides
on three circles and corners at their intersections}\label{polygons}

In this section we prove a preliminary result for classification of nets.
Roughly speaking it says that every spherical rectangle is a union
of two spherical triangles.

As we prove this result by induction, it is convenient to consider
a more general class of spherical polygons, characterized by the property
that the developing map sends their sides to three transversally intersecting
great circles and corners to the intersection points of these circles.
The net $\Gamma$ defines a triangulation of such a polygon $Q$, each face of it being
mapped by $f$ one-to-one onto one of the triangles
into which the three circles partition the sphere.
This triangulation satisfies the following properties:\newline
(P1) Each vertex inside $Q$ has degree $4$;\newline
(P2) All boundary vertices, other than corners of $Q$, have degree $3$.\newline
Combining this triangulation with its mirror copy, we obtain a triangulation $\mathbf T$ of the sphere
satisfying the following properties:\newline
(S1) $\mathbf T$ is symmetric with respect to a circle $S$ contained in the 1-skeleton of $\mathbf T$;\newline
(S2) Each vertex of $\mathbf T$ has even degree, and all its vertices not contained in $S$ have degree $4$.\newline
It is easy to show that the nets of spherical polygons with all sides mapped to three transversal circles
and all corners to intersection points of those circles are in one-to-one correspondence with triangulations
of the sphere satisfying (S1) and (S2).

Two nets are combinatorially equivalent if they can be obtained from each other by an orientation-preserving homeomorphism (mapping corners to corners and sides to sides) preserving the initial corner.

If $C$ is any of the three circles, its preimage in $Q$ is called $C$-{\em net}, denoted $\Gamma_C$.
An {\em arc} of the net $\Gamma_C$ (or an arc of $\Gamma$ if $C$ is not specified) is a connected component
of $\Gamma_C\setminus\partial Q$.
Since $f$ is a local homeomorphism on the interior of $Q$,
an arc may be homeomorphic to either an open interval with
both ends on the boundary of $Q$ (possibly, at the same corner of $Q$) or a circle in the interior of $Q$.
We'll show below (see Corollary \ref{nocircle}) that an arc of a spherical rectangle $Q$
must have at least one end at a corner of $Q$. In particular, an arc of a spherical rectangle cannot be a circle.
An arc is called {\em short} if it does not intersect other arcs of $\Gamma$.
Any two arcs of the same net $\Gamma_C$ are disjoint.

\begin{thm}\label{triangulation}
Let $Q$ be a spherical $n$-gon such that
all its sides are mapped to three transversal great circles by the
developing map, and all its corners are mapped to intersection points
of those circles. Then either $n\le 3$ or there is a triangulation of $Q$ by $n-3$ disjoint arcs of its net,
each of them connecting two non-adjacent corners of $Q$.
\end{thm}
\vspace{.1in}

{\em Proof.} It is enough to show that, unless $Q$ is a digon or a triangle, there exists an arc
of its net $\Gamma$ connecting two of its non-adjacent corners.
We prove this by induction on the number $N$ of faces of $\Gamma$.
If $\Gamma$ has one face then, since any face of $\Gamma$ is a triangle, $Q$ is a triangle.

If $N>1$ then there exists an arc $\gamma$ of $\Gamma$ adjacent to a point $p$ on the
boundary of $Q$.
Otherwise the face of $\Gamma$ adjacent to its boundary would not be simply connected.

If $\gamma$ connects two distinct corners $p$ and $q$ of $Q$ then either $p$ and $q$
are non-adjacent and we are done, or $p$ and $q$ are adjacent corners of $Q$,
and $\gamma$ partitions $Q$ into a digon and a polygon $Q'$ with the same number
of corners as $Q$ and a smaller than $N$ number of faces of its net.
By inductive hypothesis, unless $Q'$ (and thus $Q$) is a digon or a triangle,
there is an arc $\gamma'$ of $Q'$ connecting two of its non-adjacent corners.
In the latter case, $\gamma'$ is also an arc of $\Gamma$ connecting two non-adjacent corners of $Q$,
and we are done.

Suppose now that $\Gamma$ does not have any arcs connecting two corners of $Q$.
If $\gamma$ has both ends at the same point $p$ then $Q$ can be replaced by a $(n+1)$-gon $Q'$
having all sides of $Q$ plus $\gamma$ as its sides, with the number of faces of the net $\Gamma'$ of $Q'$ smaller than $N$.
There is a mapping $\iota:Q'\to Q$ such that any two distinct points of $Q'$ map to distinct points of $Q$,
except the two ends of the side $\gamma$ of $Q'$ that both map to $p$.
By the inductive hypothesis, there is an arc $\gamma'$ of $\Gamma'$
connecting two non-adjacent corners $p'$ and $q'$ of $Q'$.
Then $\iota(\gamma')$ is an arc of $\Gamma$ connecting two (possibly, adjacent) corners of $Q$,
a contradiction.

Thus we may suppose that $\gamma$ has two distinct ends $p$ and $q$ on the boundary of $Q$,
at least one of them not a corner of $Q$.
Then $\gamma$ partitions $Q$ into two polygons $Q'$ and $Q''$, with the number of corners
$n'$ and $n''$ respectively, where $n'+n''\ge n+3$.
If $n>3$ then at least one of $n'$ and $n''$ is greater than $3$.
Since both $Q'$ and $Q''$ have the number of faces of their nets smaller than $N$,
by the inductive hypothesis at least one of them has an arc $\gamma'$ of its net connecting two
non-adjacent corners. Then $\gamma'$ is also an arc of $\Gamma$ connecting two non-adjacent
corners of $Q$. This completes the proof.

\begin{figure}
\centering
\includegraphics[width=4in]{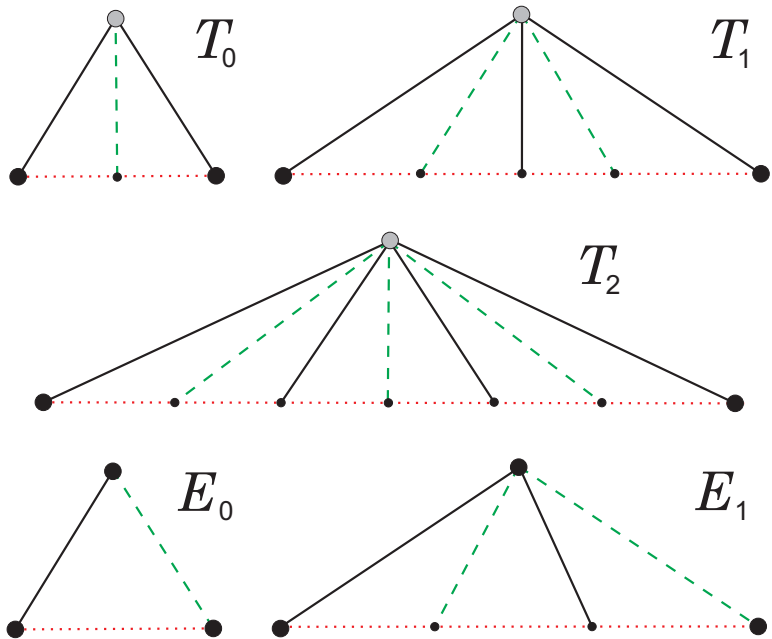}
\caption{Primitive spherical triangles $T_\mu$ and $E_\mu$.}\label{rect-triang}
\end{figure}

\begin{figure}
\centering
\includegraphics[width=5in]{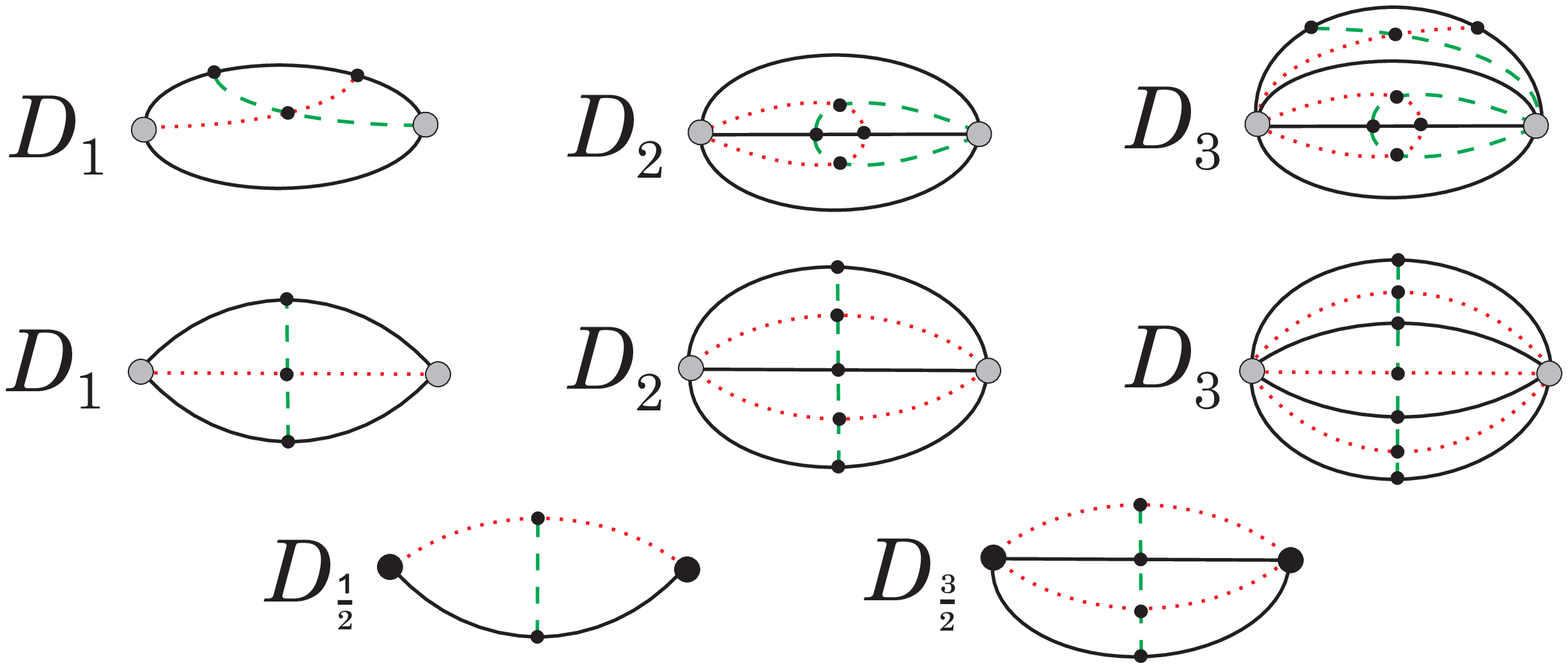}
\caption{Spherical digons.}\label{rect-digon}
\end{figure}

\begin{figure}
\centering
\includegraphics[width=5in]{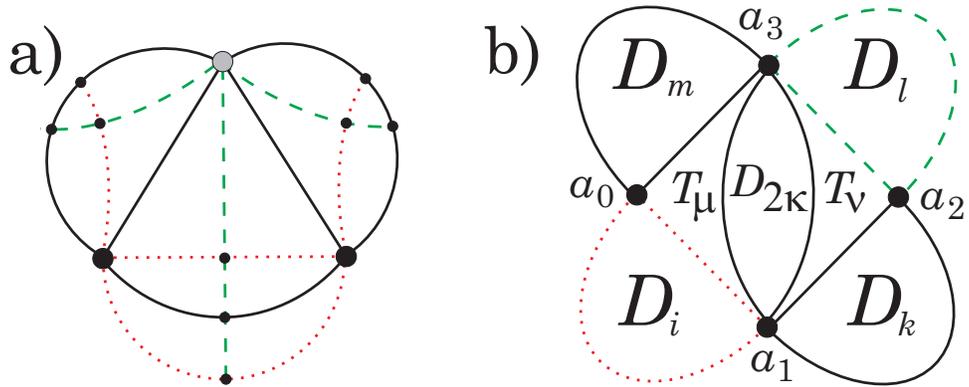}
\caption{(a) Spherical triangle $T_0$ with three digons $D_1$ attached to its sides.
(b) Nets of spherical rectangles.}\label{triang-rectang}
\end{figure}

\begin{cor}\label{nocircle}
If $Q$ is a spherical polygon satisfying the condition of Theorem \ref{triangulation}
then an arc $\gamma$ of the net of $Q$ is an open interval with at least one end at a corner of $Q$.
\end{cor}
\vspace{.1in}

{\em Proof.} It follows from classification of spherical triangles (see sections 10 and 12 of \cite{EGT2},
section 6 of \cite{EGT3} and Figs.~\ref{rect-triang}, \ref{rect-digon} and \ref{triang-rectang}a)
that an arc of the net of a spherical triangle with all sides mapped to three circles $C,\, C',\, C''$ and
all corners to intersection points of those circles must have at least one end at its corner.
Indeed, such a triangle $T$ contains a primitive triangle $T'$
(either one of triangles $T_\mu$ or one of triangles $E_\mu$, see Fig.~\ref{rect-triang})
with its apex at the intersection of
two circles, say $C$ and $C'$, and its base on the third circle $C''$.
Any arc of the net of $T'$ connects its apex with its base.
The triangle $T$ is obtained by attaching digons to the sides of $T'$ (see Fig.~\ref{triang-rectang}a).
An arc of the net of a digon $D$ either has an end at one of its corners or
the ends on both its sides, but cannot have both ends on one side of $D$.
Thus an arc of $T$ must have at least one end at its corner.
Since the intersection of $\gamma$ with any triangle $T$ of a triangulation of $Q$ by disjoint arcs of its net
connecting its non-adjacent corners is either a side of $T$ or an arc of the net of $T$, it must have at least one end at
a corner of $T$. But all corners of $T$ are also corners of $Q$.

\noindent
\section{Nets of spherical rectangles}\label{nets}

\begin{figure}
\centering
\includegraphics[width=4in]{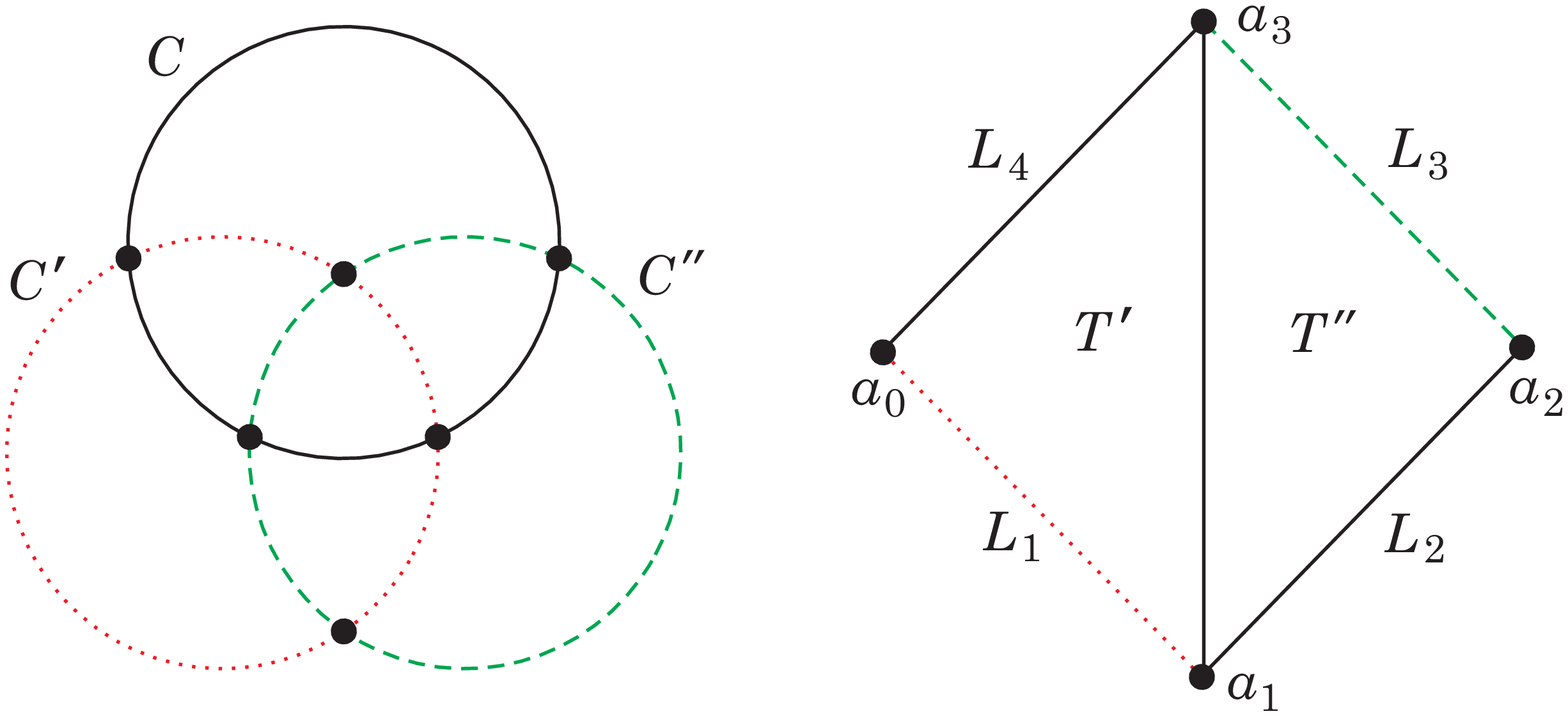}
\caption{Three circles (left) and a spherical rectangle (right).}\label{TT}
\end{figure}

As was shown in \cite{EG2} (see Proposition \ref{onecircle}) any spherical rectangle $Q$
has two opposite sides mapped to the same circle by its developing map $f$, and two other
opposite sides mapped to distinct circles.
Thus there are two types of marked spherical rectangles: in the first type
the images of $L_2$ and $L_4$ belong to the same circle,
and in the second type the images of $L_1$ and $L_3$ belong to the same circle.
It is enough to classify spherical rectangles of the first type,
as all rectangles of the second type can be obtained from those of the first type
by orientation-reversing isometry preserving the marked corner.
\vspace{.1in}

\noindent{\bf Assumption 1.} {\em Unless stated otherwise, all spherical
rectangles below are assumed to be of the first type.}
\vspace{.1in}

Let $C$ be the circle to which two sides $L_2$ and $L4$ of a spherical rectangle $Q$ are mapped,
and let $C'$ and $C''$ be the circles to which the sides $L_1$ and $L_3$ of $Q$ are mapped (see Fig.~\ref{TT}).

Theorem \ref{triangulation} implies that there is an arc of the net $\Gamma$ of $Q$
connecting two opposite corners of $Q$. Such an arc must be mapped to the circle $C$,
since two opposite corners of $Q$ are mapped to intersection points of $C$ with
two distinct circles other than $C$.
This implies that $Q$ cannot have two arcs of $\Gamma$ connecting two pairs of its
opposite corners: such arcs would have an intersection point inside $Q$,
while any two arcs mapped to the same circle $C$ are disjoint.
\vspace{.1in}

\noindent{\bf Assumption 2.} {\em Unless stated otherwise, we choose the initial corner $a_0$ of
a marked spherical rectangle $Q$ so that there is
an arc $\gamma$ of $\Gamma$ connecting the corners $a_1$ and $a_3$ of $Q$.}
\vspace{.1in}

Such an arc $\gamma$ partitions $Q$ into two spherical triangles $T'$ and $T''$,
where $T'$ has an integer corner at $a_3$ and the base $L_1$ mapped to $C'$,
while $T''$ has an integer corner at $a_1$ and the base $L_3$ mapped to $C''$ (see Fig.~\ref{TT}).
We'll show below (see Remark \ref{angles}) that the angles of such rectangle $Q$ satisfy
the inequality $A_0+A_2+2\le A_1+A_3$.

Any rectangle of the first (resp., second) type with an arc of its net
connecting its corners $a_0$ and $a_2$ can be obtained from a rectangle of the second (resp, first) type
with an arc of its net connecting its corners $a_1$ and $a_3$
by choosing $a_1$ instead of $a_0$ as an initial corner, and relabeling the corners accordingly.
The angles of such rectangle satisfy the inequality $A_1+A_3+2\le A_0+A_2$.
Thus it is enough to classify spherical rectangles satisfying Assumptions 1 and 2.

\begin{thm}\label{classification}
Let $Q$ be a marked spherical rectangle satisfying Assumptions 1 and 2.
Then $Q$ is a union of two primitive triangles $T_\mu$ and $T_\nu$ having integer angles $\mu+1$
and $\nu+1$, respectively,
digon $D_{2\kappa}$ with the sides mapped to $C$ having common sides with both $T_\mu$ and $T_\nu$,
digon $D_i$ with the sides mapped to $C'$ attached to the base of $T_\mu$,
digon $D_l$ with the sides mapped to $C''$ attached to the base of $T_\nu$,
digon $D_m$ attached to the remaining side of $T_\mu$, and digon $D_k$ attached to the remaining side of $T_\nu$.
The sides of $D_k$ and $D_m$ are mapped to $C$.

Here $\mu,\nu,\kappa,i,k,l,m$ are non-negative integers satisfying $i\mu=l\nu=0$,
that is, $i>0$ only if $\mu=0$, $l>0$ only if $\nu=0$.
The value $0$ for $i,k,l,m,\kappa$ means there is no digon attached.
\end{thm}
\vspace{.1in}

{\em Proof.}
Assumptions 1 and 2 imply that an arc $\gamma$ of the net $\Gamma$ of $Q$ connecting its corners
$a_1$ and $a_3$ partitions $Q$ into two spherical triangles $T'$ and $T''$, where $T'$ has an integer
corner with an angle $\mu+1$ at $a_3$, and $T''$ has an integer corner with an angle $\nu+1$ at $a_1$,
for some non-negative integers $\mu$ and $\nu$.

Classification of spherical triangles with one integer corner (see sections 10 and 12 of \cite{EGT2})
implies that the triangle $T'$ (resp., $T''$) is combinatorially equivalent to a
primitive triangle $T_\mu$ (resp., $T_\nu$) having an angle $\mu+1$ (resp., $\nu+1$) at its integer corner,
with digons attached to its sides (see Figs.~\ref{rect-triang}, \ref{rect-digon} and \ref{triang-rectang}a).
No digons may be attached to the base $L_1$ of $T'$ (resp., the base $L_3$ of $T''$)
if $\mu>0$ (resp., $\nu>0$).
Each digon $D_n$ has equal integer angles $n$ at its two corners.

The sides of $T_\mu$ and $T_\nu$ are mapped to $C$ and cannot contain preimages of
intersection points of $C$ with either $C'$ or $C''$, other than the corners of $Q$.
This implies that the union of digons attached
to the sides of $T_\mu$ and $T_\nu$ and having $\gamma$ as their common side is a digon $D_{2\kappa}$
with even integer angles $2\kappa$ at its two corners.

Thus a net of a spherical rectangle satisfying Assumptions 1 and 2 must have the structure
shown schematically in Fig.~\ref{triang-rectang}b. This proves Theorem \ref{classification}.
\vspace{.1in}

\begin{rmk}\label{angles}
The angles at the corners $a_0,a_1,a_2,a_3$ of a marked spherical rectangle $Q$ in Theorem \ref{classification}
have the integer parts
\begin{equation}\label{inequality}
A_0=i+m,\; A_1=i+k+\nu+1+2\kappa,\; A_2=k+l,\; A_3=l+m+\mu+1+2\kappa,
\end{equation}
respectively. In particular, $A_0+A_2+2\le A_1+A_3$.
For a marked spherical rectangle (of either first or second type) with an arc of its net
connecting its corners $a_0$ and $a_2$, the integer parts of its angles satisfy $A_1+A_3+2\le A_0+A_2$.
\end{rmk}
\vspace{.1in}

\begin{rmk}\label{shortarc}
Theorem \ref{classification} implies that a spherical rectangle $Q$ satisfying Assumptions 1 and 2
has at least one short arc $\gamma$ connecting its corners $a_1$ and $a_3$,
and that all such short arcs are mapped to the same arc $\beta$ of the circle $C$
with the ends at the intersection points of $C$ with $C'$ and $C''$.
\end{rmk}

\begin{df}\label{df:theta}
The angle $\theta\in(0,1)$ between the circles $C'$ and $C''$ is defined as $1-\alpha$
where $\alpha$ is the length of any short arc $\gamma$ of the net of
a spherical rectangle $Q$ connecting its opposite corners, divided by $\pi$.
Alternatively, $\alpha$ is the length of the arc $\beta$ of $C$ to which $\gamma$ is mapped,
divided by $\pi$.
\end{df}

{\em Proof of Proposition \ref{isometry}.}
We want to show that two marked spherical rectangles $Q$ and $Q'$
with equivalent nets $\Gamma$ and $\Gamma'$,
and the same angle $\theta$, are congruent, i.e.,
there is an orientation-preserving
isometry $Q\to Q'$ mapping the initial corner
$a_0$ of $Q$ to the initial corner $a'_0$ of $Q'$.
According to Definition \ref{df:theta},
the developing map $f:Q\to\bC$ maps the sides of $Q$
to three great circles $C$, $C'$, $C''$,
such that both $C'$ and $C''$ are orthogonal to $C$,
so that the sides $L_2$ and $L_4$ of $Q$ are mapped to $C$,
the sides $L_1$ and $L_3$ of $Q$ are mapped to $C'$ and $C''$,
respectively, and any short arc
of $\Gamma$ connecting the corners $a_1$ and $a_3$ of $Q$
is mapped to an arc $\beta$ of $C$ of length
$\pi\alpha$ where $\alpha=1-\theta$,
the ends $P$ and $R$ of $\beta$ being the images of $a_1$ and $a_3$,
respectively. In particular, $P$ and $R$ are intersection
points of $C$ with $C'$ and $C''$, respectively.
Similarly, the developing map $g:Q'\to\bC$ maps the sides of
$Q'$ to three great circles $S$, $S'$, $S''$,
such that both $S'$ and $S''$ are orthogonal to $S$,
and any short arc $\gamma'$ of $\Gamma'$ connecting
the corners $a'_1$ and $a'_3$ of $Q'$ is mapped to an arc $\beta'$
of $S$ of length $\pi\alpha$, the ends $P'$ and $R'$
of $\beta'$ being the images of $a'_1$ and $a'_3$, respectively.
Applying if necessary rotation of $\bC$,
we may assume that $S=C$, $S'=C'$, $S''=C''$, $\beta'=\beta$,
$P'=P$, and $R'=R$.
Then the equivalence of the nets $\Gamma$ and $\Gamma'$
implies that the developing maps $f$ and $g$
send the corresponding faces, edges and vertices
of partitions of $Q$ and $Q'$ to the same triangles, segments
and intersection points of the partition of $\bC$ by
the three circles $C$, $C'$ and $C''$.
In particular, $a_0$ and $a'_0$ are mapped by $f$
and $g$ to the same point of $\bC$.

\begin{thm}\label{existence}
For any non-negative integers $\mu,\nu,\kappa,i,k,l,m$
satisfying $i\mu=l\nu=0$ there is a unique,
up to combinatorial equivalence, net $\Gamma$ of the type
described in Theorem \ref{classification}.
For any such net $\Gamma$ and any $\theta\in(0,1)$
there exists a unique spherical rectangle $Q=Q(\Gamma,\theta)$
having $\Gamma$ as its net, sides mapped to three circles
$C$, $C'$, $C''$, and a short arc of length $\pi(1-\theta)$
connecting its corners $a_1$ and $a_3$.
\end{thm}
\vspace{.1in}

{\em Proof.} To define the net $\Gamma$, we start with a digon $D_{2\kappa}$
obtained by combining $\kappa$ copies of digon $D_2$ shown
in the middle of the first row
of Fig.~\ref{rect-digon}. If $\kappa=0$ then there is
no such digon, and we proceed with
gluing together triangles $T_\mu$ and $T_\nu$
(see Fig.~\ref{rect-triang}) so that the side
of each of these triangles
that follows its base in the counterclockwise
cyclic order becomes their common side.
If $\kappa>0$ then triangles $T_\mu$ and $T_\nu$
are attached to opposite sides of $D_{2\kappa}$
so that the side of each of these triangles that follows
its base in the counterclockwise cyclic
order becomes its common side with $D_{2\kappa}$.
The integer corner of $T_\mu$ (resp., $T_\nu$) coincides
with a non-integer corner of $T_\nu$ (resp., $T_\mu$).

Next, we attach digons $D_k$ and $D_m$, obtained by combining
$k$ and $m$ copies, respectively,
of the digon $D_1$ shown in the left side of the
first row of Fig.~\ref{rect-digon}, so that any
two adjacent digons have either a common short side
or a common long side, and so that each of
the two resulting digons has at least one short side,
as is shown in examples of $D_2$ and $D_3$ in
the first row of Fig.~\ref{rect-digon}.
Then digons $D_k$ and $D_m$ are attached
to the free sides of
$T_\nu$ and $T_\mu$, respectively.
The free sides of these triangles are preceding
their respective
bases in the counterclockwise cyclic order.
If $k=0$ (resp., $m=0$) then no digon $D_k$
(resp., $D_m$) is attached.

Finally, if $\mu=0$ and $i>0$ (resp., $\nu=0$ and $l>0$)
then a digon $D_i$ (resp., $D_l$),
obtained by combining $i$ and $l$ copies, respectively,
of the digon $D_1$ shown in the second row
of Fig.~\ref{rect-digon}, is attached to the base
of $T_\mu$ (resp., $T_\nu$).
Examples of such digons are shown in the
second row of Fig.~\ref{rect-digon}.


If we label the sides of $T_\mu$ and $T_\nu$ by $C$,
the base of $T_\mu$ by $C'$ and the base of $T_\nu$
by $C''$ then all edges of $\Gamma$ can be uniquely
labeled so that the sides of each of its triangles
are labeled by three distinct labels.
Consider the standard sphere $\bC$ with three great circles $C$, $C'$
and $C''$ such that $C$ is orthogonal to both $C'$ and $C''$,
and one of the two complementary angles
between $C'$ and $C''$ is $\theta$,
the other one being $\alpha=1-\theta$.
Then there is a mapping of $Q$ to $\bC$,
locally one-to-one everywhere except at the corners,
which is unique up to a homeomorphism of $Q$
preserving all vertices and edges of $\Gamma$ and a rotation
of the sphere preserving the three circles,
such that any arc of the boundary of $D_{2\kappa}$
(or a common side of $T_\mu$ and $T_\nu$ if $\kappa=0$)
maps to an arc of $C$ of length $\pi\alpha$,
and each edge of $\Gamma$ maps to an arc of the circle
corresponding to its label. This defines on $Q$ a metric
of the spherical rectangle $Q(\Gamma,\theta)$.
\vspace{.1in}

\begin{figure}
\centering
\includegraphics[width=5in]{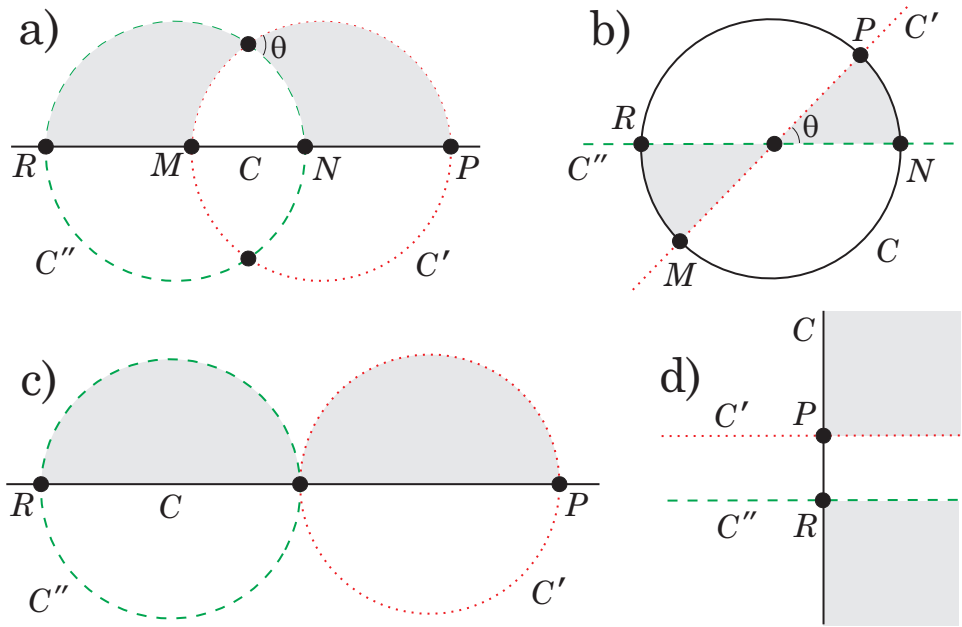}
\caption{Deformation of the three-circle configuration}\label{rect-3circles}
\end{figure}

\begin{figure}
\centering
\includegraphics[width=5in]{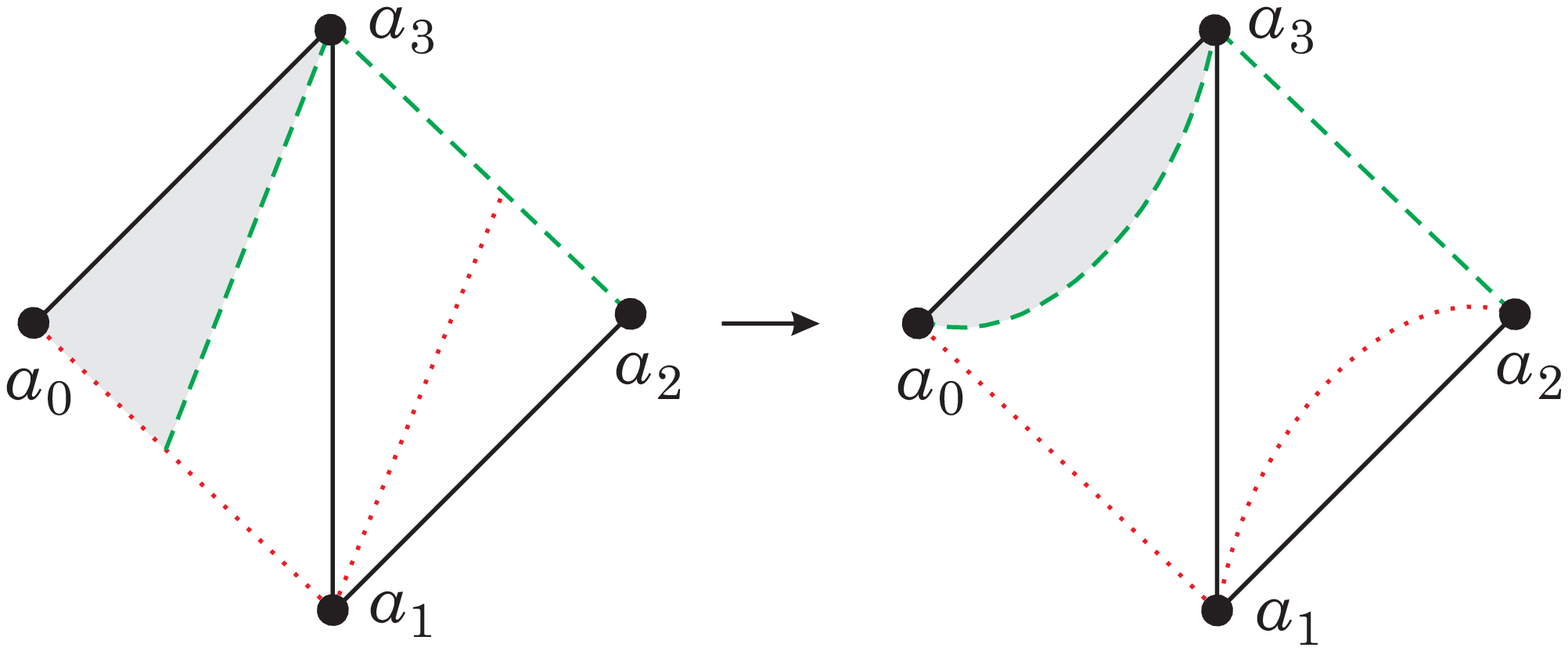}
\caption{Deformation of spherical rectangles with the angles $\frac12,\,\frac32,\,\frac12,\,\frac32$}\label{rectangle1313}
\end{figure}

\begin{figure}
\centering
\includegraphics[width=4in]{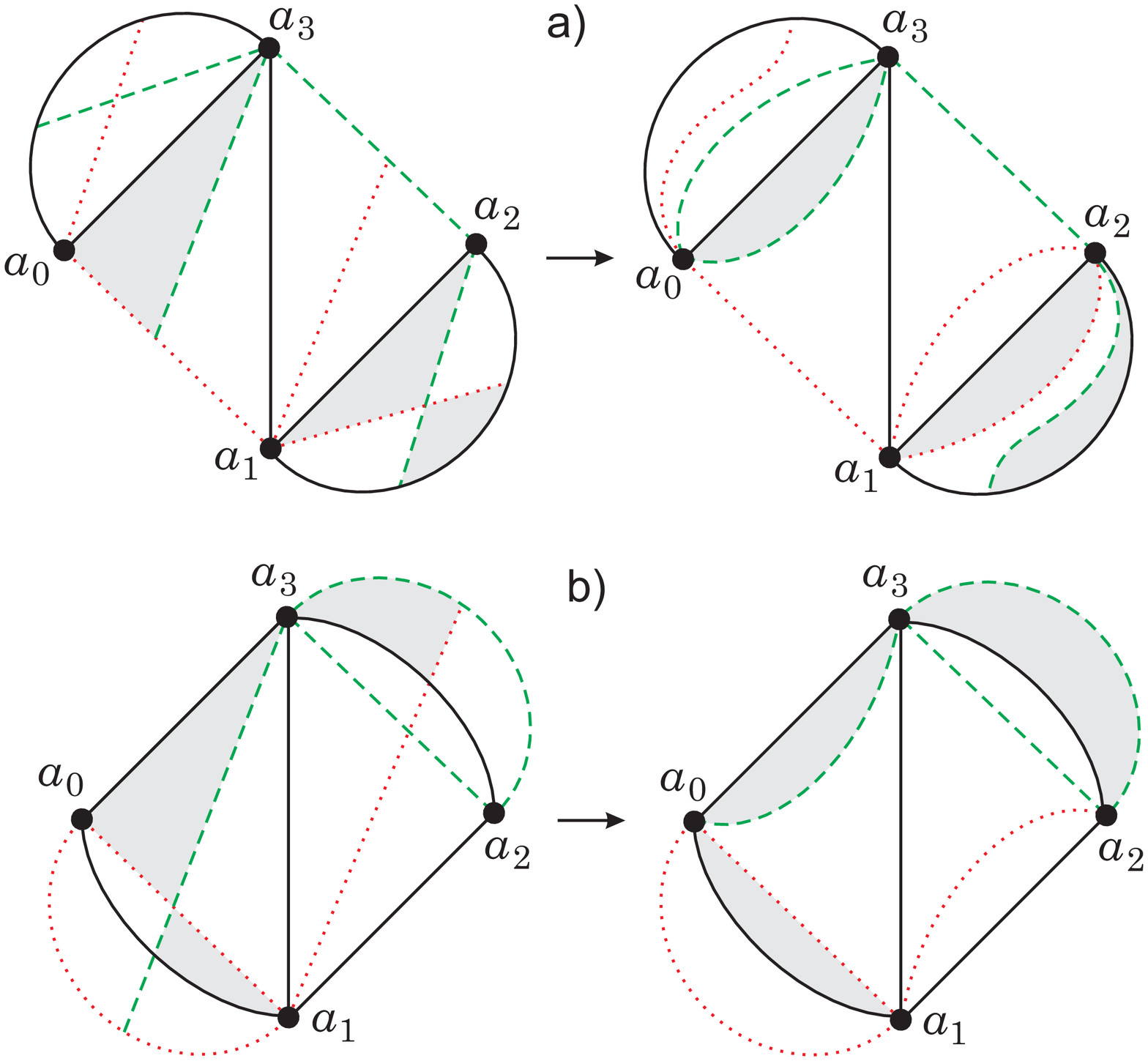}
\caption{Deformation of spherical rectangles with the angles $\frac32,\,\frac52,\,\frac32,\,\frac52$}\label{rectangle3535}
\end{figure}

\begin{figure}
\centering
\includegraphics[width=4in]{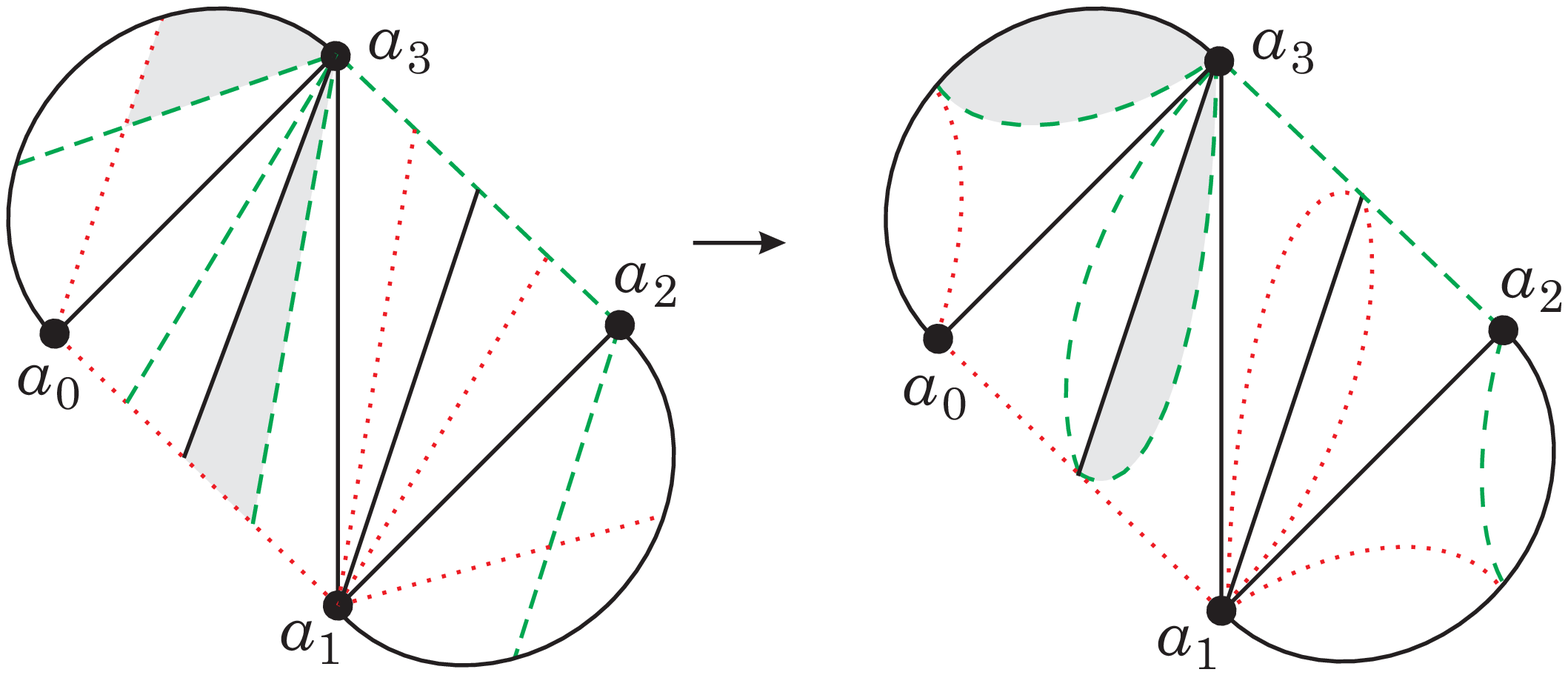}
\caption{Deformation of spherical rectangles with the angles $\frac32,\,\frac72,\,\frac32,\,\frac72$}\label{rectangle3737}
\end{figure}

\begin{figure}
\centering
\includegraphics[width=5in]{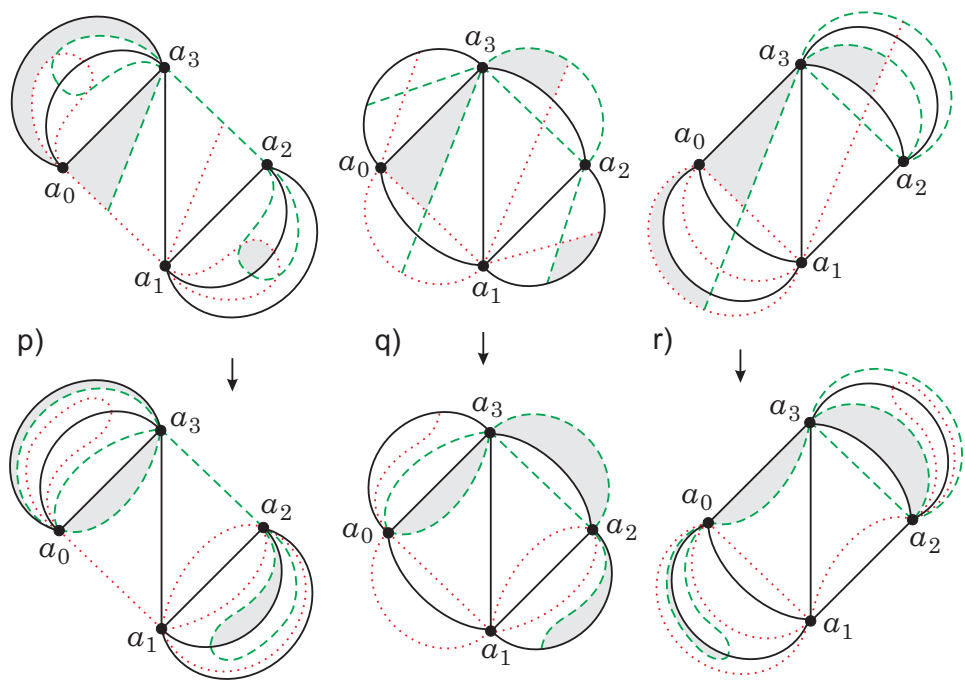}
\caption{Three combinatorially distinct nets of spherical rectangles with the angles $\frac52,\,\frac72,\,\frac52,\,\frac72$, and their deformations}\label{rectangle5757}
\end{figure}

\noindent
\section{Limits at $\theta=0$ and $\theta=1$ of spherical rectangles with the given net}\label{limits}

Each admissible set of integers $\mu,\nu,\kappa,i,k,l,m$ in Theorem \ref{classification}
defines a net $\Gamma$ and the corresponding one-parametric family of marked spherical
rectangles $Q(\Gamma,\theta)$ satisfying Assumptions 1 and 2,
parameterized by the angle $\theta$ between the circles $C'$ and $C''$ (see Definition \ref{df:theta}).
These two circles intersect the circle $C$ at the right angle.
In Figs.~\ref{rect-3circles}a and \ref{rect-3circles}b, two projections of the three
circles are shown, and two of the triangles of the partition of the sphere defined by these
circles are shaded.

Let $P$ and $R$ be the images of the corners $a_1$ and $a_3$, respectively,
of a marked spherical rectangle $Q$, so that the arc $\beta=PR$ of $C$ is the image of
a short arc $\gamma$ of the net $\Gamma$ of $Q$ connecting $a_1$ and $a_3$.
Then the shaded areas in Figs.~\ref{rect-3circles}a and \ref{rect-3circles}b
contract to arcs when $\theta\to 0$ and expand to half-disks when $\theta\to 1$.
If all circles remain geodesic, then $C'$ and $C''$ converge to the same circle when
$\theta\to 0$ and when $\theta\to 1$. However, applying a linear-fractional transformation
depending on $\theta$, so that the arc $MN$ of $C$ in Figs.~\ref{rect-3circles}a and \ref{rect-3circles}b
contracts to a point while the arc $PR$ does not, we can obtain in the limit $\theta\to 1$ a non-geodesic configuration
shown in two different projections in Figs.~\ref{rect-3circles}c and \ref{rect-3circles}d,
where the shaded areas are the limits of the shaded areas in Figs.~\ref{rect-3circles}a and \ref{rect-3circles}b.

\begin{example}\label{1313} {\rm
A net of a spherical rectangle $Q$ with the angles
\newline
$(1/2,3/2,1/2,3/2)$
considered in \cite{EG2} is shown in Fig.~\ref{rectangle1313} (left).
The shaded area corresponds to preimage of the shaded areas in Figs.~\ref{rect-3circles}a and \ref{rect-3circles}b.
The arc connecting $a_1$ and $a_3$ is mapped to the arc $PR$ of $C$, the corners $a_0$ and $a_2$ are mapped to $M$ and $N$, respectively.
When $\theta\to 0$, the sides $L_2$ and $L_4$ of $Q$ are contracted to points, while the distance between them
has a positive limit.
Thus the modulus of $Q$ has limit $0$ as $\theta\to 0$ (see \cite{EGT2}, Section 15, and \cite{EG2}).
When $\theta\to 1$, applying a linear-fractional deformation depending on $\theta$ to the three-circle
configuration, and the corresponding transformation to $Q$
(this does not change the modulus of $Q$ which is a conformal invariant) we can get in the limit
a non-geodesic circular rectangle shown in Fig.~\ref{rectangle1313} (right).
Thus the modulus of $Q$ tends to a finite positive value when $\theta\to 1$.
Computation in \cite{EG2} shows that this value is $K\approx 0.630963$.}
\end{example}

\begin{thm}\label{limit}
Let $\Gamma$ be one of the nets described in Theorem \ref{classification}, and $\theta\in(0,1)$.
Then the modulus of $Q(\Gamma,\theta)$ tends to $0$ when $\theta\to 0$, and to a finite positive value when $\theta\to 1$.
\end{thm}
\vspace{.1in}

{\em Proof.}
When $\theta\to 0$, none of the arcs
connecting $a_1$ and $a_3$ contract, while the triangle
$T_\mu$ (resp., $T_\nu$) contains either a short arc
of $\Gamma$ or a side $L_4$ (resp., $L_2$) of $Q$
that maps to either the arc $MR$ or the arc $NP$ of $C$,
connecting its apex $a_3$ (resp., $a_1$) with a point $p\ne a_1$
(resp., $p\ne a_3$) on its base.
These two arcs contract to points when $\theta\to 0$.
Hence the distance between the sides $L_1$ and $L_3$ of $Q$ tends to $0$ when $\theta\to 0$.
At the same time, there are no short arcs of $\Gamma$ having one end on $L_2$ and another end on $L_4$,
other than those connecting $a_1$ and $a_3$ which do not contract as $\theta\to 0$.
Thus the distance between the sides $L_2$ and $L_4$ does not tend to $0$ as $\theta\to 0$.
This implies that the modulus of $Q$ tends to $0$ as $\theta\to 0$ (see \cite{EGT2}, Section 15).

When $\theta\to 1$, the short arcs connecting $a_1$ and $a_3$ contract,
thus both the distance between
$L_1$ and $L_3$ and the distance between $L_2$ and $L_4$ tend to $0$.
To understand the limit of the modulus of $Q$,
we apply a linear-fractional transformation depending on $\theta$
as in Example \ref{1313} to the three-circle configuration,
so that the short arc $MN$ of $C$
contracts, while the short arc $PR$ does not.
In the non-geodesic limit (see Fig.~\ref{rect-3circles}c) the circles $C'$ and $C''$ become tangent (when their tangency point is mapped to $\infty$
as in Fig.~\ref{rect-3circles}d, they become parallel lines).
All short arcs of $\Gamma$ connecting $a_1$ and $a_3$ map to $PR$, and all arcs of $C$ connecting the apex of a triangle $T_\mu$ (resp., $T_\nu$)
with a point on its base, map to either $MR$ or $NP$.
Since neither $MR$ nor $NP$ contracts when $\theta\to 1$, and $PR$ does not contract after the linear-fractional transformation,
the distances between opposite sides of $Q$ do not tend to $0$ in the limit.
Thus $Q$ converges to a non-geodesic circular rectangle (see Figs.~\ref{rectangle1313}, \ref{rectangle3535}, \ref{rectangle3737},
and \ref{rectangle5757}) and the modulus of $Q$ tends to a
finite positive value.
\vspace{.2in}

\noindent
{\bf Remarks on computation of limit moduli $K$.}
The boundary of the limit rectangle described in the proof
of Theorem~\ref{limit} is mapped by developing map into three straight lines
(see Fig.~\ref{rect-3circles}d). This allows to represent the
developing map by the Schwarz--Christoffel formula. Condition that the
points $P$ and $Q$ are on the same vertical line imposes one real
equation which permits to determine the modulus of the rectangle
$K$. See \cite{EG2} where the simplest example is described in all detail.
The number of solutions to this equation is the number of nets
with given angles. To determine which solution corresponds to which net,
we use the evident inequalities between the moduli of degenerate
rectangles (shown in the right of Figs.~\ref{rectangle3535}, \ref{rectangle3737}, and in the bottom
of Fig.~\ref{rectangle5757}).

\begin{figure}
\centering
\includegraphics[width=4in]{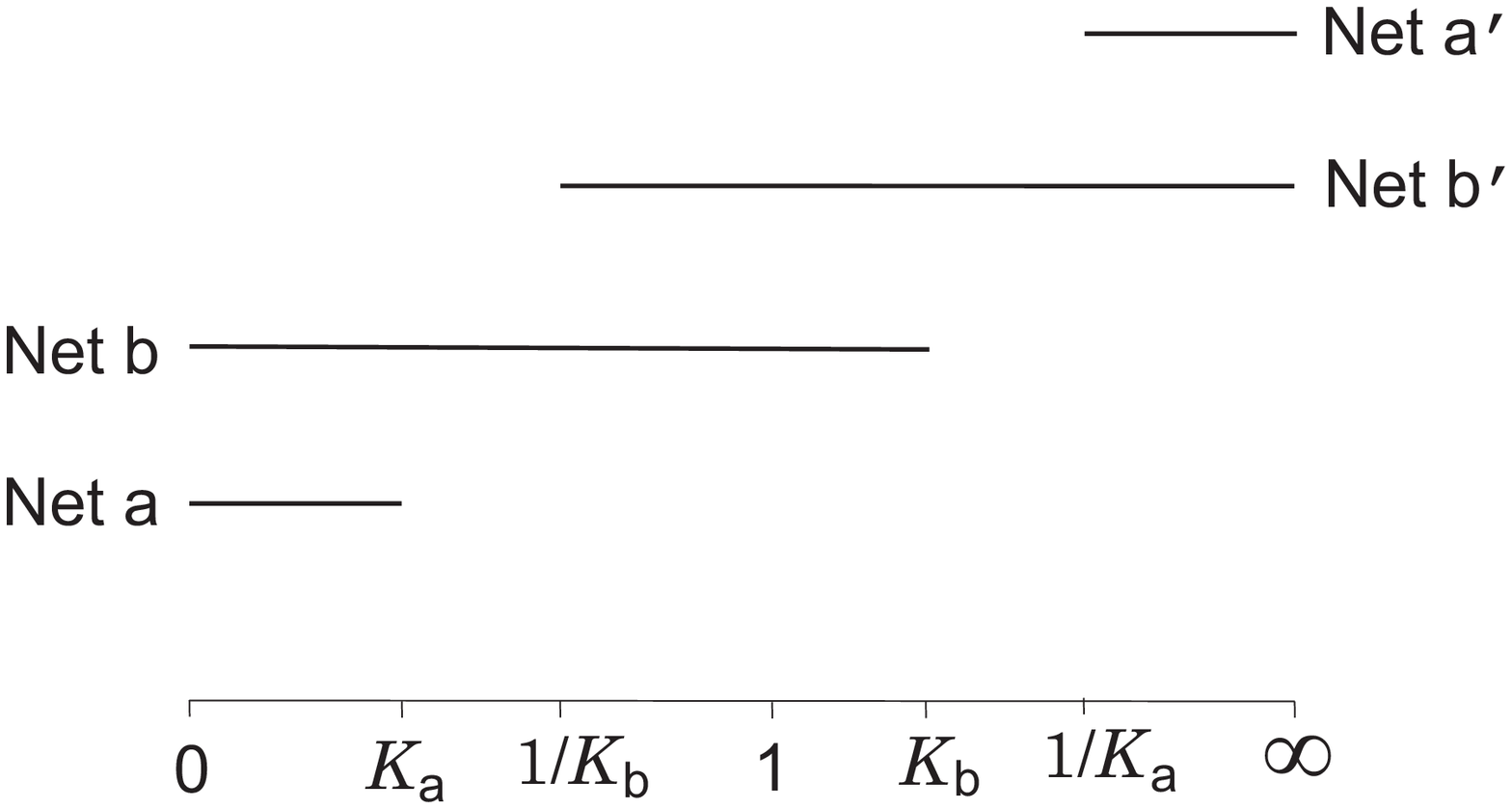}
\caption{Existence of spherical rectangles with the angles $\frac32,\,\frac52,\,\frac32,\,\frac52$}\label{KK}
\end{figure}

\begin{example}\label{3535} {\rm
Two combinatorially distinct nets of spherical rectangles with the angles
$(3/2,5/2,3/2,5/2)$,
and the nets of their non-geodesic limits when $\theta\to 1$,
are shown in Figs.~\ref{rectangle3535}a and \ref{rectangle3535}b.
The moduli $K_a$ and $K_b$ of the limiting rectangles are
$K_a\approx 0.5433144$ and $K_b\approx 1.193606$, respectively.
Fig.~\ref{KK} shows schematically the areas
of existence of these spherical rectangles
(Nets a and b) and their involution-symmetric rectangles
(Nets $\mathrm{a}'$ and $\mathrm{b}'$) for different values
of the modulus $K$.}
\end{example}

\begin{figure}
\centering
\includegraphics[width=5in]{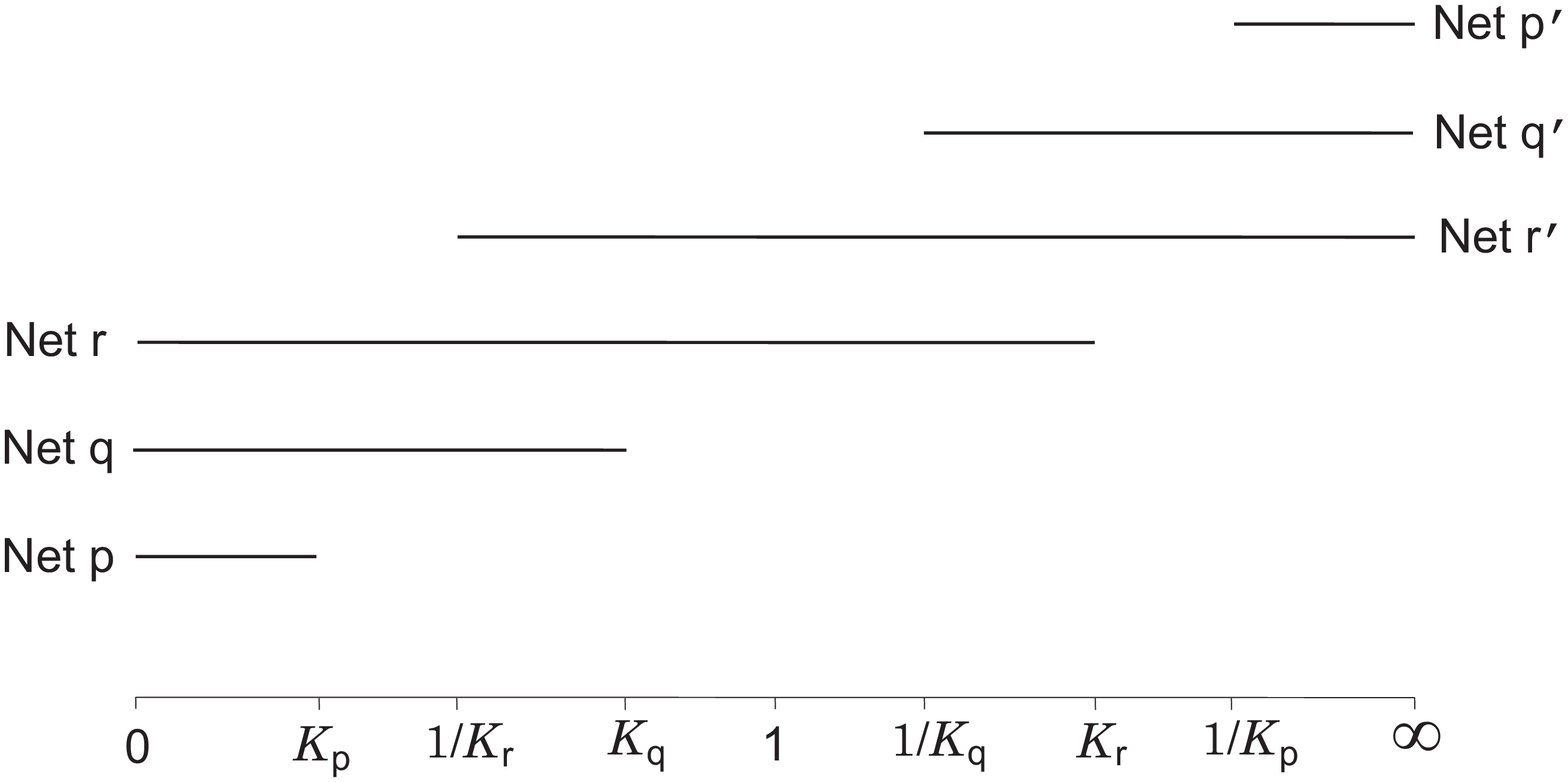}
\caption{Existence of spherical rectangles with the angles $\frac52,\,\frac72,\,\frac52,\,\frac72$}\label{KKK}
\end{figure}

\begin{example}\label{5757} {\rm
Three combinatorially distinct nets of spherical
rectangles with the angles $(5/2,7/2,5/2,7/2)$,
and the nets of their non-geodesic limits when $\theta\to 1$,
are shown in Figs.~\ref{rectangle5757}p, \ref{rectangle5757}q,
\ref{rectangle5757}r. The moduli $K_p$, $K_q$, $K_r$ of the
limiting rectangles are $K_p\approx 0.476966$,
$K_q\approx 0.887943$, $K_r\approx 1.458956$,
respectively.
Fig.~\ref{KKK} shows schematically the areas of existence
of these spherical rectangles
(Nets p, q, r) and their involution-symmetric rectangles
(Nets $\mathrm{p}'$, $\mathrm{q}'$, $\mathrm{r}'$) for different values
of the modulus $K$.}
\end{example}

\begin{example}\label{3737} {\rm
A net of a spherical rectangle with the angles\newline $(3/2,7/2,3/2,7/2)$,
and the net of its non-geodesic limit when $\theta\to 1$, is shown in Fig.~\ref{rectangle3737}.
The modulus $K$ of the limiting rectangle is $K\approx 0.4173$.}
\end{example}

\noindent
{\bf Remarks and conjectures.} Each of the three nets p, q, r in Fig.~\ref{KKK}
produces a continuous family where the modulus $K$ can be
arbitrarily small. So for
sufficiently small $K$ there are at least three different
marked spherical quadrilaterals with the angles $(5/2,7/2,5/2,7/2)$.
Similarly, one can conclude from Fig.~\ref{KKK} that
there are at least two quadrilaterals with modulus $K\in(K_{\mathrm{p}},
1/K_{\mathrm{r}}$ at least three for $K\in (1/K_{\mathrm{r}},K_{\mathrm{q}})$,
at least two for $K\in(K_{\mathrm{q}},1/K_{\mathrm{q}})$,
at least three for $K\in(1/K_{\mathrm{q}},K_{\mathrm{r}})$,
at least two for $K\in(K){\mathrm{r}},1/K_{\mathrm{p}})$,
and at least three for $K>1/K_{\mathrm{p}}$.

Similar conclusions apply to Fig.~\ref{KK}. Our computations show that in
fact these lower estimates are equalities, in all cases which we computed.
Actually $K$ is a monotone function of $\theta$ in all these cases, but
we do not expect this monotonicity to hold for all angles.

So Theorem \ref{theorem2} gives only
lower estimates for the number of quadrilaterals with given angles
and modulus, when this modulus is small or large. This lower estimate
is $N$ or $M_j$ (see (\ref{K1}), (\ref{K2}), (\ref{N}).
We conjecture that these lower estimates are exact, and that we have
equality for large and small moduli.
All computed examples confirm this.

\section{Counting nets of spherical rectangles with given angles}

We want to answer the following question: Given four non-negative integers $A_0,\dots,A_3$,
how many nets of marked spherical rectangles with the integer parts $A_0,\dots,A_3$ of the angles
at their corners $a_0,\dots,a_3$ do exist?

It is enough to answer this question for the spherical rectangles of the first type (satisfying Assumption 1),
since all rectangles of the second type can be obtained then by an involution preserving the marked corner.
Also, we may assume that $A_0+A_2+2\le A_1+A_3$, which is true for the marked spherical rectangles
satisfying Assumption 2 (see Remark \ref{angles}).
The answer for the rectangles with $A_1+A_3+2\le A_0+A_2$, having an arc connecting corners $a_0$ and $a_2$,
can be obtained then by a different choice of the initial corner.

We start with listing operations on the nets which do not change the angles of a spherical rectangle $Q$.
Assuming notation of Theorem \ref{classification} and Remark \ref{angles}, we have expressions (\ref{inequality})
for the angles of $Q$.

\medskip
\noindent{\bf Operation I.} If $\kappa>0,\; i\ge 2,\; l=0,\; \mu=0$ then
$$\kappa\mapsto\kappa-1,\; \nu\mapsto\nu+4,\; i\mapsto i-2,\; m\mapsto m+2.$$
The inverse operation is possible when $l=0,\; m\ge 2,\; \mu=0,\; \nu\ge 4$.

\smallskip
\noindent{\bf Operation II.} If $\kappa>0,\; i=0,\; l\ge 2,\; \nu=0$ then
$$\kappa\mapsto\kappa-1,\; \mu\mapsto\mu+4,\; k\mapsto k+2,\; l\mapsto l-2.$$
The inverse operation is possible when $i=0,\; k\ge 2,\; \mu\ge 4,\; \nu=0$.

\smallskip
\noindent{\bf Operation III.} If $\kappa>0,\; i=1,\; l=0,\; \mu=0$ then
$$\kappa\mapsto\kappa-1,\; \mu\mapsto 1,\; \nu\mapsto\nu+3,\; i\mapsto 0,\; m\mapsto m+1.$$
The inverse operation is possible when $i=l=0,\; m\ge 1,\; \mu=1,\; \nu\ge 3$.

\smallskip
\noindent{\bf Operation IV.} If $\kappa>0,\; i=0,\; l=1,\; \nu=0$ then
$$\kappa\mapsto\kappa-1,\; \mu\mapsto\mu+3,\; \nu\mapsto 1,\; k\mapsto k+1,\; l\mapsto 0.$$
The inverse operation is possible when $i=l=0,\; k\ge 1,\; \mu\ge 3,\; \nu=1$.

\smallskip
\noindent{\bf Operation V.} If $\kappa>0,\; i=l=0$ then
$$\kappa\mapsto\kappa-1,\; \mu\mapsto\mu+2,\; \nu\mapsto\nu+2.$$
The inverse operation is possible when $i=l=0,\; \mu\ge 2,\; \nu\ge 2.$

\smallskip
\noindent{\bf Operation VI.} If $i>0,\; l>0,\; \mu=\nu=0$ then
$$i\mapsto i-1,\; k\mapsto k+1,\; l\mapsto l-1,\; m\mapsto m+1.$$
The inverse operation is possible when $k>0,\; m>0,\; \mu=\nu=0.$

\begin{lemma}\label{kappa}
Let $Q$ be a marked spherical rectangle satisfying Assumptions 1 and 2.

a) If $\kappa>0$ and $\min(i,l)=0$ then there is a unique operation among Operations I-V
that is applicable to $Q$ and results in a rectangle with the same angles as $Q$,
with $\kappa$ reduced by $1$.

b) If $\min(i,l)=0$ then at most one operation among inverses to Operations I-V
is applicable to $Q$. If $\min(i,l)>0$ then neither Operations I-V nor their inverses
are applicable to $Q$.

c) If $\min(i,l)>0$, then iteration of Operation VI applied to $Q$
results in a rectangle with the same angles as $Q$, same $\kappa,\,\mu,\,\nu$, and $\min(i,l)=0$.
\end{lemma}
\vspace{.1in}

{\em Proof.}
We start with a proof of (a), assuming $\kappa>0$.
If $i=l=0$ then Operation V is applicable.
If $i=1,\; l=0$ then $\mu=0$ and Operation III is applicable.
If $i\ge 2,\; l=0$ then $\mu=0$ and Operation I is applicable.
If $i=0,\; l=1$ then $\nu=0$ and Operation IV is applicable.
If $i=0,\; l\ge 2$ then $\nu=0$ and Operation II is applicable.
It is easy to check that only one of the operations I-V is applicable in each of these cases.

To prove (b) note first that inverses to Operations I-V are only possible when $\min(i,l)=0$.
Next, for given $\mu$ and $\nu$, conditions on $\mu$ and $\nu$ for applicability of the operations inverse to Operations I-V may
hold for at most one of these operations.

If $\min(i,l)>0$ then Operation VI reduces $i,\, l$ and $\min(i,l)$ by 1,
and does not change $\kappa,\,\mu,\,\nu$, which proves (c).

\begin{cor}\label{exist}
For given $A_0,\dots,A_3$, the set of values of $\kappa$ that may appear in the nets of marked
spherical rectangles of the first type with $\delta\ge 1$ is either empty (in which case spherical
rectangles with such angles do not exist) or an interval $[0,\kappa_{\rm max}]$, for some integer $\kappa_{\rm max}\ge 0$
depending on $A_0,\dots,A_3$.
In the latter case, there are exactly $\kappa_{\rm max}+1$ combinatorially distinct
nets of marked spherical rectangles of the first type with given $A_0,\dots,A_3$ and $\min(i,l)=0$.
\end{cor}

For a marked spherical rectangle $Q$ satisfying Assumption 2, define
\begin{equation}\label{delta1}
\delta=(A_1+A_3-A_0-A_2)/2=2\kappa+1+(\mu+\nu)/2.
\end{equation}
Then $\delta\ge 1$ is either integer or half-integer.

\begin{lemma}\label{kappamax}
Let $Q$ be a marked spherical rectangle with $\min(i,l)=0$ satisfying Assumptions 1 and 2.
The net of $Q$ cannot be obtained from a net of some other spherical rectangle
by one of the operations inverse to Operations I-V if and only if
one of the following twelve conditions is satisfied:\newline
{\rm (a)} $\mu=\nu=0$;\newline
{\rm (b)} $\mu=0,\;\nu=1$;\newline
{\rm (c)} $\mu=1,\;\nu=0$;\newline
{\rm (d)} $\mu=\nu=1$;\newline
{\rm (e)} $\mu=0,\;\nu=2$;\newline
{\rm (f)} $\mu=2,\;\nu=0$;\newline
{\rm (g)} $\mu=0,\;\nu=3$;\newline
{\rm (h)} $\mu=3,\;\nu=0$;\newline
{\rm (i)} $\mu=1,\;\nu\ge 3,\;m=0$;\newline
{\rm (j)} $\mu\ge 3,\;\nu=1,\;k=0$;\newline
{\rm (k)} $\mu=0,\;\nu\ge 4,\; m\le 1$;\newline
{\rm (l)} $\mu\ge 4,\;\nu=0,\; k\le 1$.\newline
For given $A_0,\dots,A_3$ with $A_1+A_3\ge A_0+A_2+2$, a net of a marked spherical rectangle with $\min(i,l)=0$
satisfying Assumptions 1 and 2 may satisfy at most one of these conditions.\newline
At most one (up to combinatorial equivalence) net of a marked spherical rectangle $Q$
satisfying Assumptions 1 and 2, with given $A_0,\dots,A_3$ and $\min(i,l)=0$, may satisfy any of these conditions.
The value of $\kappa$ for such rectangle $Q$ is
\begin{equation}\label{eq:kappamax}
\left[\min\left(\frac{A_1-1}2,\frac{A_3-1}2,\frac{\delta-1}2\right)\right].
\end{equation}
\end{lemma}
\vspace{.1in}

{\em Proof.}
One can easily check case by case that none of the operations inverse to Operations I-V can be applied
if and only if one of the conditions (a)-(l) is satisfied. Note that it is enough to assume $\mu\le\nu$,
and to check that with this assumption none of the operations inverse to Operations I,III,V can be applied
if and only if one of the conditions (a),(b),(d),(e),(g),(i),(k) is satisfied.
The case $\mu\ge\nu$ follows by rotation of the net exchanging $a_0$ with $a_2$, $a_1$ with $a_3$,
$\mu$ with $\nu$, $i$ with $l$, and $k$ with $m$.

Let now $Q$ be a marked spherical rectangle satisfying Assumptions 1 and 2, with given $A_0,\dots,A_3$,
$A_1+A_3\ge A_0+A_2+2$, $\min(i,l)=0$, and $\mu\le\nu$.

If $Q$ satisfies (a) then $\delta=2\kappa+1$ is an odd integer, $A_1\ge\delta$, $A_3\ge\delta$.
If, in addition, $l=0$ then $A_2=k$ and $A_3=\delta+m$, thus $i=A_0-m=A_0-A_3+\delta=(A_0+A_1-A_2-A_3)/2$, and
the net of $Q$ is completely determined by its angles. Such a net exists when $i\ge 0$ thus $A_0+A_1\ge A_2+A_3$.
The case when $i=0$ is treated similarly, and a net with $i=0$ satisfying (a) exists when $A_0+A_1\le A_2+A_3$.
The net with $i=l=0$ satisfying (a) exists when $A_0+A_1=A_2+A_3$.

If $Q$ satisfies (b) then $\delta=2\kappa+\frac32$ thus $2\delta\equiv 3 \!\!\mod 4$,
$A_1\ge\delta+\frac12$, $A_3\ge\delta-\frac12$, $l=0$, $A_2=k$, $A_3=m+\delta-\frac12$ thus $m=A_3-\delta+\frac12$,
$A_0=i+m=i+A_3-\delta+\frac12$ thus $i=A_0-A_3+\delta-\frac12=(A_0+A_1-A_2-A_3-1)/2$.
A net satisfying (b) is completely determined by its angles. It exists when $i\ge 0$, thus $A_0+A_1\ge A_2+A_3+1$.
Similarly, if $Q$ satisfies (c) then $\delta=2\kappa+\frac32$, and its net is completely determined by its angles.
It exists when $A_1\ge\delta-\frac12$, $A_3\ge\delta+\frac12$, $A_2+A_3\ge A_0+A_1+1$.

If $Q$ satisfies (d) then $\delta=A_3-A_0=A_1-A_2=2\kappa+2$ is a positive even integer, $i=l=0$.
Then $A_0=m$, $A_2=k$, and the net is completely determined by its angles. It exists when $A_0+A_1=A_2+A_3$.

If $Q$ satisfies (e) then $\delta=2\kappa+2$ is a positive even integer,
$A_1\ge\delta+1$, $A_3\ge\delta-1$, $l=0$, $A_0=i+m$,
$A_2=k$, $A_1=i+k+\delta+1$ thus $i=A_1-A_2-\delta-1=(A_0+A_1-A_2-A_3-2)/2$,
$m=A_3-\delta+1$, and the net is completely determined by its angles.
It exists when $i\ge 0$ thus $A_0+A_1\ge A_2+A_3+2$.
Similarly, if $Q$ satisfies (f) then $\delta=2\kappa+2$, and its net is completely determined by its angles.
It exists when $A_1\ge\delta-1$, $A_3\ge\delta+1$, $A_2+A_3\ge A_0+A_1+2$.

If $Q$ satisfies (g) then $\delta=2\kappa+\frac52$ thus $2\delta\equiv 1 \!\!\mod 4$,
$A_1\ge\delta+\frac32$, $A_3\ge\delta-\frac32$, $l=0$, $A_2=k$, $A_3=m+\delta-\frac32$ thus $m=A_3-\delta+\frac32$,
$A_0=i+m=i+A_3-\delta+\frac32$ thus $i=A_0-A_3+\delta-\frac32=(A_0+A_1-A_2-A_3-3)/2$.
A net satisfying (g) is completely determined by its angles. It exists when $i\ge 0$, thus $A_0+A_1\ge A_2+A_3+3$.
Similarly, if $Q$ satisfies (h) then $\delta=2\kappa+\frac52$, and its net is completely determined by its angles.
It exists when $A_1\ge\delta+\frac32$, $A_3\ge\delta-\frac32$, $A_2+A_3\ge A_0+A_1+3$.

If $Q$ satisfies either (i), (j), (k) or (l) then $\delta=2\kappa+3$ is an odd integer, $\delta\ge 3$.
If $Q$ satisfies (i) then $A_3=\delta-1$, $i=l=0$, $A_0=0$, $A_2=k$.
If $Q$ satisfies (j) then $A_1=\delta-1$, $i=l=0$, $A_0=m$, $A_2=0$.
If $Q$ satisfies (k) then $A_3=\delta-2+m$, $l=0$, $A_0=i+m$, $A_2=k$.
A net satisfying (k) exists when $i=A_0-m=A_0-A_3+\delta-2=(A_0+A_1-A_2-A_3-4)/2\ge 0$, thus $A_0+A_1\ge A_2+A_3+4$.
If $Q$ satisfies (l) then $A_1=\delta-2+k$, $i=0$, $A_0=m$, $A_2=k+l$.
A net satisfying (l) exists when $l=A_2-k=A_2-A_1+\delta-2=(A_2+A_3-A_0-A_1-4)/2\ge 0$, thus $A_2+A_3\ge A_0+A_1+4$.
A net satisfying either (i), (j), (k) or (l) is completely determined by its angles.

\begin{lemma}\label{kappamaxplus1}
A marked spherical rectangle $Q$ satisfying Assumptions 1 and 2, with given $A_0,\dots,A_3$ and $\min(i,l)=0$
exists if and only if
\begin{equation}\label{eq:exist}
\min\left(A_1,A_3,\delta\right)\ge 1,
\end{equation}
 and there are exactly
\begin{equation}\label{eq:kappamaxplus1}
\left[\min\left(\frac{A_1+1}2,\frac{A_3+1}2,\frac{\delta+1}2\right)\right]
\end{equation}
combinatorially distinct nets of such spherical rectangles.
\end{lemma}
\vspace{.1in}

{\em Proof.}
It follows from Corollary \ref{exist} and Lemma \ref{kappamax} that existence of a rectangle $Q$ satisfying conditions of Lemma \ref{kappamaxplus1}
implies that the number in (\ref{eq:kappamax}) is non-negative, which is equivalent to (\ref{eq:exist}),
and that the number of combinatorially distinct nets of such rectangles in that case is the number in (\ref{eq:kappamax}) plus $1$,
which is the number in (\ref{eq:kappamaxplus1}). Thus we have only to prove that for any $A_0,\dots,A_3$ satisfying (\ref{eq:exist})
there exists a rectangle $Q$ satisfying conditions of Lemma \ref{kappamaxplus1}.

It follows from Corollary \ref{exist} that, if a rectangle with given $A_0,\dots,A_3$ satisfying conditions of Lemma \ref{kappamaxplus1}
exists, there exists also a rectangle with $\kappa=0$ with the same angles satisfying the same conditions.
We are going to construct a net of such a rectangle for any $A_0,\dots,A_3$ satisfying (\ref{eq:exist}).

There are three possible cases: (i) $A_1>A_2$ and $A_3>A_0$;
(ii) $A_1>A_2$ and $A_3\le A_0$; (iii) $A_1\le A_2$ and $A_3>A_0$.
Note that $A_1\le A_2$ and $A_3\le A_0$ is not possible because $\delta=(A_1+A_3-A_0-A_2)/2\ge 1$.\newline
In case (i) let $i=l=0$, $m=A_0$, $k=A_2$, $\mu=A_3-m-1=A_3-A_0-1$, $\nu=A_1-k-1=A_1-A_2-1$.\newline
In case (ii) let $\mu=l=0$, $m=A_3-1$, $k=A_2$, $i=A_0-m=A_0-A_3+1$, $\nu=A_1-i-k-1=A_1+A_3-A_0-A_2-2=2(\delta-1)$.\newline
In case (iii) let $\nu=i=0$, $k=A_1-1$, $m=A_0$, $l=A_2-k=A_2-A_1+1$, $\mu=A_3-l-m-1=A_1+A_3-A_0-A_2-2=2(\delta-1)$.\newline
This completes the proof of Lemma \ref{kappamaxplus1}.

\begin{df}\label{df:special}
{\rm A marked spherical rectangle is {\em special} if $\delta$ is an odd integer
and either $A_1\ge\delta>0$ and $A_3\ge\delta>0$ or $A_0\ge -\delta>0$ and $A_2\ge -\delta>0$.}
\end{df}

\begin{lemma}\label{special}
A marked spherical rectangle $Q$ satisfying Assumptions 1 and 2 with given $A_0,\dots,A_3$ is
special if an only if there exists a rectangle with the same angles whose net has $\mu=\nu=0$.
For given $A_0,\dots,A_3$ satisfying conditions of Definition \ref{df:special} with $\delta>0$, there are
\begin{equation}\label{extra-special}
\min(A_0,A_1-\delta,A_2,A_3-\delta)
\end{equation}
special rectangles satisfying Assumptions 1 and 2, with $\mu=\nu=0$ and $\min(i,l)>0$.
\end{lemma}
\vspace{.1in}

{\em Proof.}
If $\min(i,l)>0$ for the net of $Q$ then $\mu=\nu=0$.
From Lemma \ref{kappa} (c), iteration of Operation VI applied to $Q$
results in a unique rectangle $Q_0$ with the same angles as $Q$, $\mu=\nu=0$ and $\min(i,l)=0$.
Thus $Q_0$ satisfies Lemma \ref{kappamax} (a).
If $\min(i,l)=0$ for the net of $Q$ then, from Corollary \ref{exist}, there is a unique
rectangle with the same angles as $Q$ satisfying one of the conditions (a)-(l) of Lemma \ref{kappamax}.
One can easily check that $\delta$ is an odd integer, $A_1\ge\delta$ and $A_3\ge\delta$ only in case (a)
of Lemma \ref{kappamax}.
Conversely, a rectangle satisfying condition (a) of Lemma \ref{kappamax} is clearly special,
thus any rectangle $Q$ with the same angles is also special.
Finally, the number in (\ref{extra-special}) is obtained by counting distinct rectangles
that can be obtained from a rectangle with $\mu=\nu=0$ by iterating Operation VI and its inverse
(compare with Lemma 11.2 in \cite{EGT2}).

\begin{thm}\label{count-all}
Let $Q$ be a marked spherical rectangle satisfying Assumptions 1 and 2, with given $A_0,\dots,A_3$.
If $Q$ is not special then there are (\ref{eq:kappamaxplus1}) combinatorially distinct nets of marked
spherical rectangles satisfying Assumptions 1 and 2 with the same angles as $Q$.
If $Q$ is special then the number in (\ref{eq:kappamaxplus1}) is $(1+\delta)/2$,
and there are additionally (\ref{extra-special}) nets
of marked spherical rectangles satisfying Assumptions 1 and 2 with the same angles as $Q$,
thus the total number of combinatorially distinct nets is
\begin{equation}\label{total-special}
\min(A_0+\frac{1+\delta}2,A_1+\frac{1-\delta}2,A_2+\frac{1+\delta}2,A_3+\frac{1-\delta}2).
\end{equation}
\end{thm}
\vspace{.1in}

{\em Proof.} This follows from Lemmas \ref{kappamaxplus1} and \ref{special}.
\vspace{.1in}

{\em Proof of Theorems \ref{theorem1} and \ref{theorem2}.}
Lemmas \ref{kappamaxplus1} and \ref{special} imply that a marked spherical
rectangle satisfying Assumptions 1 and 2 exists if and only if (\ref{eq:exist}) holds.
Since this condition is symmetric with respect to $A_1$ and $A_3$, it remains true
for spherical rectangles of the second type satisfying $A_0+A_2<A_1+A_3$,
as any such rectangle can be obtained from a rectangle satisfying Assumptions 1 and 2
by a reflection preserving $a_0$ and $a_2$, exchanging $a_1$ and $a_3$.

If $Q$ is a marked spherical rectangle of either first of second type satisfying $A_0+A_2> A_1+A_3$,
replacing $a_0$ by $a_1$ as an initial corner, and relabeling the corners accordingly,
results in a marked spherical rectangle $Q'$ of either second or first type, with the integer parts of the angles
$(A'_0,A'_1,A'_2,A'_3)=(A_1,A_2,A_3,A_0)$ and $\delta'=(A'_1+A'_3-A'_0-A'_2)/2=-\delta$.
Applying the above arguments to $Q'$ we see that
a marked spherical rectangle $Q$ with $A_0+A_2> A_1+A_3$ exists if and only if
\begin{equation}\label{eq:exist1}
\min\left(A_0,A_2,-\delta\right)\ge 1,
\end{equation}
Combining (\ref{eq:exist}) and (\ref{eq:exist1}) we get the statement of Theorem \ref{theorem1}.
The statement of Theorem \ref{theorem2} follows from (\ref{eq:kappamaxplus1}) and (\ref{total-special})
applied to either $Q$ or $Q'$ in a similar way.

\begin{lemma}\label{symmetric}
Let $Q$ be a marked spherical rectangle with the angles satisfying $A_0=A_2,\;A_1=A_3$.
Then there is an orientation-preserving isometry $\rho: Q\to Q$ such that $\rho(a_0)=a_2$ and $\rho(a_1)=a_3$.
\end{lemma}
\vspace{.1in}

{\em Proof.} Due to Proposition \ref{isometry}, it is enough to prove that the net $\Gamma$
of $Q$ is symmetric with respect to a transformation exchanging $a_0$ with $a_2$ and $a_1$ with $a_3$,
and that this symmetry of $\Gamma$ maps any short arc of $\Gamma$ connecting $a_1$ with $a_3$ is mapped to a
(possibly, different) short arc connecting $a_1$ with $a_3$.
To show this, we have only to check (assuming that $Q$ satisfies Assumptions 1 and 2)
that $\mu=\nu$, $i=l$, and $k=m$ in the notations of Theorem \ref{classification}.

Suppose first that $\mu=\nu=0$. Then $A_0=A_2$ implies $i+m=k+l$, and $A_1=A_3$ implies $i+k=m+l$
(see (\ref{inequality}) in Remark \ref{angles}).
Adding up these two equalities yields $i=l$, and subtracting them yields $k=m$.

If $\mu>0$ and $\nu>0$ then $i=l=0$, thus $A_0=m$, $A_2=k$, $A_1=k+1+2\kappa+\nu$, and
$A_3=m+1+2\kappa+\mu$.
Since $A_0=A_2$, we have $k=m$, then $A_1=A_3$ yields $\mu=\nu$.

If $\mu>0$ but $\nu=0$ then $i=0$, thus $A_0=m$, $A_1=k+1+2\kappa$, $A_2=k+l$, and
$A_3=l+m+1+2\kappa+\mu$. Since $A_0=A_2$, we have $m=k+l$, thus
$A_3=k+2l+1+2\kappa+\mu>A_1$, a contradiction.
Similarly, $\mu=0$ and $\nu>0$ is not possible.
This completes the proof.
\vspace{.1in}

\begin{thm}\label{thm:symmetric}
A spherical rectangle $Q$ with the angles at two of its opposite corners equal $\alpha$,
and the angles at two other opposite corners equal $\beta$, exists if and only if
$|\beta-\alpha|\ge 1$. If $\beta-\alpha$ is even then there are $|\beta-\alpha|$
non-isometric spherical rectangles with these angles, $|\beta-\alpha|/2$ of them
satisfying Assumption 1. If $\beta-\alpha$ is odd
then there are $\alpha+\beta$ non-isometric spherical rectangles with these angles,
$(\alpha+\beta)/2$ of them satisfying Assumption 1.
\end{thm}

This follows from Lemma \ref{symmetric} and Theorems \ref{theorem1} and \ref{theorem2}.

\begin{example}{\rm
The net in Fig.~\ref{rectangle1313} is special, with $A_0=A_2=0$, $A_1=A_3=1$, $\delta=1$.
According to (\ref{total-special}) and Theorem \ref{thm:symmetric}, there is a unique net of a marked spherical rectangle of the first type with these angles.

The two nets in Fig.~\ref{rectangle3535} are special with $A_0=A_2=1$, $A_1=A_3=2$, $\delta=1$.
According to (\ref{total-special}) and Theorem \ref{thm:symmetric}, there are two nets of marked spherical rectangles of the first type with these angles.

The net in Fig.~\ref{rectangle3737} is not special, with $A_0=A_2=1$, $A_1=A_3=3$, $\delta=2$.
According to (\ref{eq:kappamaxplus1}) and Theorem \ref{thm:symmetric}, there is a unique net of a marked spherical rectangle of the first type with these angles.

The three nets in Fig.~\ref{rectangle5757} are special with $A_0=A_2=2$, $A_1=A_3=3$, $\delta=1$.
According to (\ref{total-special}) and Theorem \ref{thm:symmetric},
there are three nets of marked spherical rectangles of the first
type with these angles.}
\end{example}

{\em Department of Mathematics, Purdue University,
West Lafayette Indiana, 47907 USA}

\begin{thebibliography}{1}
\bibitem{khov} A. Khovanskii,On a class of systems of transcendental equations, Soviet Math. Dokl., 22 (1980), 762--765.
\bibitem{Lin2}  Ching-Li Chai, Chang-Shou Lin and Chin-Lung Wang,
Mean field equations, hyperelliptic curves and modular forms: I,
arXiv:1502.03297.
\bibitem{Darboux} G. Darboux, Sur une \'equation lin\'eaire,
C. R. Acad. Sci. Paris, t. 94 (1882) 1645--1648.
\bibitem{E} A. Eremenko,
Metrics of positive curvature with conic singularities on
the sphere, Proc. Amer. Math. Soc. 132 (2004), 3349--3355
\bibitem{EG2}
 A. Eremenko and A. Gabrielov, On metrics of curvature 1 with four conic
singularities on tori and on the sphere, arXiv:1508.06510.
\bibitem{EGSV} A. Eremenko, A. Gabrielov, M. Shapiro and A. Vainshtein,
Rational functions and real Schubert calculus, Proc. AMS, 134 (2006), no. 4,
949--957.
\bibitem{EGT1} A. Eremenko, A. Gabrielov and V. Tarasov,
Metrics with conic singularities and spherical polygons,
arXiv:1405.1738
\bibitem{EGT2} A. Eremenko, A. Gabrielov and V. Tarasov,
Metrics with four conic singularities and spherical quadrilaterals,
arXiv:1409.1529
\bibitem{EGT3}
 A. Eremenko, A. Gabrielov and V. Tarasov, Spherical quadrilaterals
with three non-integer angles, arXiv:1504.02928.
\bibitem{Finch} S. Finch, Mathematical constants, Cambridge UP, Cambridge, 2003.
\bibitem{Hermite} Ch. Hermite, Sur l'\'equation de Lam\'e,
Extrait de feuilles authographi\'ees du Course d'Analyse de l'\'Ecole
Polytechnique, $1^{\mathrm{re}}$ Division, 1872-73, $32^e$ le\c{c}on.
Oeuvres, t. III, p. 118--122. Paris, Gauthier-Villars, 1912.
\bibitem{Klein2} F. Klein, Vorlesungen \"uber das Insider
und die Auf\"snug der Clonking Von f\"often Grade, Birkh\"auser Verlag,
Basel, 1993.
\bibitem{Lin1} Chang-Shou Lin and Chin-Lung Wang,
Elliptic functions Green functions and the mean field equations on
tori, Ann. Math., 172 (2010) 911--954.
\bibitem{LT} F. Luo and G. Tian, Liouville equation and spherical
convex polytopes,
Proc. Amer. Math. Soc. 116 (1992), no. 4, 1119.1129.
\bibitem{MP} G. Mondello and D. Panov, Spherical metrics
with conical singularities on a 2-sphere: angle constraints,
arXiv:1505.01994.
\bibitem{Schneps} L. Schneps, Dessins d'enfants on the Riemann sphere, in:
The Grothendieck theory of dessins d'enfants (Luminy, 1993), 47--77,
Cambridge Univ. Press, Cambridge, 1994.
\bibitem{sp} P. Speissegger, The Pfaffian closure of an o-minimal structure,
J. f\"ur die reine und angew. Math., 508 (1999) 189--211.
\bibitem{Troy2} M. Troyanov, Prescribing curvature on compact
surfaces with conical
singularities, Trans. Amer. Math. Soc., 324 (1991) 793--821.
\bibitem{vdd} L. van den Dries, Tame topology and o-minimal structures, London Math. Soc. Lecture Notes Series 248,
Cambridge University Press, Cambridge 1998.
\bibitem{VV} E. Van Vleck, A determination of the number of real and imaginary
roots of the hypergeometric series, Trans. Amer. Math. Soc., 3 (1902)
110--131.
\end{thebibliography}
\end{document}